\documentclass[11pt]{article}
\usepackage{latexsym}
\usepackage[all]{xy}
\usepackage{amsmath}
\usepackage{amssymb}
\usepackage{amsfonts}

\newtheorem{thm}{Th\'eor\`eme}[section]
\newtheorem{prop}[thm]{Proposition}
\newtheorem{lem}[thm]{Lemme}

\newtheorem{df}[thm]{D\'efinition}
\newtheorem{cor}[thm]{Corollaire}

\newtheorem{rmk}[thm]{Remarque}

\begin{document}

\title{\textbf{Alg\'ebrisation des vari\'et\'es analytiques complexes et cat\'egories d\'eriv\'ees}}
\bigskip
\bigskip

\author{\bigskip\\
Bertrand To\"en et Michel Vaqui\'e\\
\small{Institut de Math\'ematiques de Toulouse}\\
\small{UMR CNRS 5219} \\
\small{Universit\'{e} Paul Sabatier, Bat 1R2}
\small{Toulouse Cedex 9}}

\date{Mars 2007}

\maketitle

\begin{abstract}
Soit $X$ un espace analytique complexe compact et lisse. Nous d\'emontrons que
$X$ est alg\'ebrisable si et seulement si sa dg-cat\'egorie d\'eriv\'ee
coh\'erente born\'ee est satur\'ee. 
\end{abstract}

\medskip

\tableofcontents

\bigskip

\section{Introduction}

Pour un sch\'ema $X$ propre et lisse sur $Spec\, \mathbb{C}$ (et plus g\'en\'eralement
pour un espace alg\'ebrique complexe propre et lisse $X$), la cat\'egorie 
triangul\'ee $D_{parf}(X)$ des complexes parfaits de $\mathcal{O}_{X}$-modules
poss\`ede de remarquables propri\'et\'es de finitudes. Dans \cite{bv} il est d\'emontr\'e
que $D_{parf}(X)$ est \emph{satur\'e}, au sens ou ses Ext globaux sont
de dimension finie et de plus tout foncteur cohomologique de type fini
$H : D_{parf}(X)^{op} \longrightarrow Vect$ est repr\'esentable. De plus, 
Bondal et Van den Bergh d\'emontrent aussi que si $X$ est une surface
complexe compacte qui ne poss\`ede pas des courbes compactes 
alors $D_{parf}(X)$ n'est pas satur\'ee. 
Ils sugg\`erent alors que toute vari\'et\'e complexe compacte $X$ dont la
cat\'egorie d\'eriv\'ee parfaite est satur\'ee est alg\'ebrisable (voir \cite[Rem . 5.6.2]{bv}). 
Le but de ce travail est de d\'emontrer cette assertion lorsque l'on remplace
la cadre \'etriqu\'e des cat\'egories triangul\'ees par celui plus flexible des dg-cat\'egories.

Pour $X$ une vari\'et\'e complexe compacte, la cat\'egorie triangul\'ee
$D_{parf}(X)$ est la cat\'egorie homotopique d'une dg-cat\'egorie naturelle
$L_{parf}(X)$ (voir \cite{tova} et la d\'efinition \ref{ddg}). De plus, une notion 
de \emph{dg-cat\'egorie satur\'ee} est introduite dans \cite{tova}, 
qui est proche de celle de cat\'egorie triangul\'ee satur\'ee (voir l'appendice
A pour une comparaison). Le th\'eor\`eme principal de ce travail est le suivant.

\begin{thm}\label{ti}
Soit $X$ un espace analytique compact et lisse (i.e. une
vari\'et\'e complexe compacte et lisse). Alors $X$ est alg\'ebrisable
(par un espace alg\'ebrique) si et seulement si la dg-cat\'egorie $L_{parf}(X)$
est satur\'ee.
\end{thm}

La n\'essit\'e de la condition est bien connue, nous rappellerons une esquisse de preuve
dans l'appendice B. La partie difficile est la suffisance, 
dont la preuve est une application du r\'esultat principal 
de \cite{tova} qui affirme le caract\`ere alg\'ebrique du champ des modules 
des objets dans une dg-cat\'egorie satur\'ee. Plus pr\'ecis\`ement, 
nous commencerons par consid\'erer l'espace
alg\'ebrique $M$ des objets simples dans la dg-cat\'egorie $L_{parf}(X)$. L'existence
de cet espace alg\'ebrique est assur\'ee par les r\'esultats de \cite{tova}
(voir le th\'eor\`eme \ref{t1}). Nous montrerons ensuite que 
l'analytifi\'e $M^{an}$ de cet espace alg\'ebrique est un espace
de modules (au sens analytique) pour les complexes parfaits et simples sur $X$
(voir Prop. \ref{p10}). Nous consid\'erons
alors le morphisme  $j : X \longrightarrow M^{an}$ qui \`a un poit $x$ de $X$ associe
le faisceau gratte-ciel $k(x)$, et nous montrons que $j$ est une immersion ouverte
(voir Prop. \ref{p11}). Enfin, en retreignant les fonctions m\'eromorphes de $M^{an}$
\`a $X$ nous d\'eduisons que $X$ est un espace de Moishezon et donc est un 
espace alg\'ebrique (voir Prop. \ref{p12}). \\

Ce travail est d\'ecoup\'e en quatres sections et trois appendices. 
Dans une premi\`ere partie nous rappelons quelques faits sur
les foncteurs de changement de sites pour les faisceaux. Nous
introduirons aussi la notion de contextes g\'eom\'etriques et de
faisceaux g\'eom\'etriques, qui nous permettront de construire
le foncteur d'analytification. La section suivante rappelle quelques
r\'esultats de la th\'eorie des cat\'egories d\'eriv\'ees en g\'eom\'etrie
analytique. Ces deux premi\`eres sections consistent essentiellement
en des rappels de r\'esultats connus, et le lecteur pourra commencer sa lacture
\`a la section \S 4 et consulter les sections \S 2,3 lorsque cela est n\'ecessaire.
La section \S 4 pr\'esente d'une part le r\'esultat d'existence d'un 
espace alg\'ebrique de modules pour les dg-modules simples sur une dg-alg\`ebre propre et lisse, 
anisi qu'une description de l'analytifi\'e de cet espace. Enfin, la section 
\S 5 pr\'esente la preuve du th\'eor\`eme principal. 

Dans les appendices, le lecteur trouvera tout d'abord une comparaison entre
les notions de cat\'egories triangul\'ees satur\'ees et de dg-cat\'egorie
satur\'ee, ce qui permet de remettre notre th\'eor\`eme dans le cadre
de \cite{bv}. L'appendice B pr\'esente tr\`es bri\`evement une preuve
de la n\'ecessit\'e du th\'eor\`eme \ref{ti}, qui est probablement un fait 
d\'ej\`a connu. Enfin, l'appendice C pr\'esente un formalisme qui permet 
de manipuler des cat\'egories d\'eriv\'ees non born\'ees \`a l'aide
de techniques de cat\'egories de mod\`eles. L'existence de cette structure
de mod\`eles sera utilis\'ee tout au long de ce travail.

\section{Rappels sur le foncteur d'analytification}

Dans cette section nous rappellerons les constructions de foncteur
de changements de sites pour les faisceaux. Nous introduirons aussi 
la notion de contextes g\'eom\'etriques et de faisceaux
g\'eom\'etriques relatifs \`a ces contextes. Nous montrerons alors
que sous certaines conditions les foncteurs de changement de sites
pr\'eservent les objets g\'eom\'etriques. Nous appliquerons cela aux
contextes alg\'ebrique et analytique, et nous obtiendrons ainsi 
un foncteur d'analytification des espaces
alg\'ebriques vers les espaces analytiques. Mises \`a part les notions
de contextes et de faisceaux g\'eom\'etriques les r\'esultats 
de cette section sont bien connus. 

\subsection{Changements de sites pour les faisceaux}

Soit $C$ une cat\'egorie poss\'edant des limites finies.  
Nous notons $Pr(C)$ la cat\'egorie des pr\'efaisceaux sur $C$, 
c'est-\`a-dire la cat\'egorie des foncteurs de $C^{op}$ \`a valeurs dans la 
cat\'egorie $Ens$ des ensembles. 

Pour tout $x\in C$ nous d\'efinissons le pr\'efaisceau 
$h_x$ par
$${h}_x(y) = C(y,x).$$
Le foncteur $h : C \longrightarrow Pr(C)$ est un 
pleinement fid\`ele et nous appelons \emph{pr\'efaisceau
repr\'esentable} tout pr\'efaisceau dans son image essentielle.  
Tout pr\'efaisceau $F\in Pr(C)$ admet une r\'esolution  par 
des pr\'efaisceaux repr\'esentables, c'est-\`a-dire peut \^etre obtenu comme 
colimite de pr\'efaisceaux repr\'e\-sentables. 

Un foncteur $f:C\to D$ d\'efinit une adjonction entre les cat\'egories
de pr\'efaisceaux 
$$\xymatrix{Pr(C) \ar@<2pt>[r]^{f_!} & Pr(D) \ar@<2pt>[l]^{f^*}} \ , $$
o\`u le foncteur adjoint \`a droite $f^*$ est d\'efini par
$f^*(F)(x) := F\bigl ( f(x)\bigr )$ pour tout pr\'efaisceau $F$ dans $SPr(D)$
et pour tout $x$ dans $C$. Le foncteur adjoint \`a gauche $f_!$ est d\'efini
sur les pr\'efaisceaux repr\'esentables par
$f_!({h}_x) := {h}_{f(x)}$ pour tout $x$ dans $C$,
et est \'etendu au pr\'efaisceaux par extension de Kan \`a gauche.

Nous rappelons qu'un foncteur $f:C \longrightarrow D$ est dit \emph{exact 
\`a gauche} s'il pr\'eserve les limites finies, c'est-\`a-dire s'il pr\'eserve 
les produits fibr\'es et l'objet final. 

\begin{prop}\label{p1}
Soient $C$ et $D$ deux cat\'egories poss\'edant des limites finies, et  $f:C\longrightarrow D$ 
un foncteur exact \`a gauche. Alors le foncteur induit
$f_! : Pr(C) \longrightarrow Pr(D)$ est 
exact \`a gauche.  
\end{prop}

\textit{Preuve:} 
Comme $C$ poss\`ede un objet final $*$ et que $f(*)\simeq *$, nous avons
$*=h_{*}$ et  $f_{!}(*)\simeq h_{f(*)}\simeq h_{*}$. 

Il nous reste \`a  montrer que un diagrammes de pr\'efaisceaux  
$\xymatrix{F\ar[r] &  H & \ar[l] G}$, le morphisme naturel
$$f_!\bigl ( F\times _H G \bigr ) \longrightarrow  
f_!(F) \times  _{f_!(H)} f_!(G) $$
est un isomorphisme dans $Pr(D)$. 
Comme les colimites d'ensembles  sont universelles (i.e. stables
par changement de bases) elles
sont aussi universelles dans 
la cat\'egorie des pr\'efaisceaux. De plus,  comme tout pr\'efaisceau
peut \^etre obtenu comme colimite de pr\'efaisceaux
repr\'esentables  nous pouvons supposer que les
pr\'efaisceaux $F$ et $G$ sont de la forme respectivement
${h}_x$ et ${h}_y$. 
Nous sommes donc ramen\'es \`a montrer que le morphisme naturel
$$f_!\bigl ({h}_x \times  _H {h}_y \bigr ) \longrightarrow 
{h}_{f(x)} \times ^h _{f_!(H)} {h}_{f(y)}$$
est un isomorphisme. 

Nous rappelons que si $x$, $y$ et $z$ sont des objets de $C$ nous
avons des isomorphismes naturels
$${h}_{x\times _z y} \simeq {h}_x \times _{{h}_z} {h}_{y}.$$
Par cons\'equent, comme le foncteur $f$ commute aux produits fibr\'es, nous trouvons le
r\'esultat dans le cas o\`u  $H$ est lui aussi repr\'esentable:
$$f_!\bigl ({h}_x \times _{{h}_z} {h}_y \bigr ) \simeq 
{h}_{f(x\times _z y)} = {h}_{f(x)\times _{f(z)} f(y)} \simeq   
f_!({h}_x) \times _{f_!({h}_z)} f_!({h}_y)\ .$$
De ceci, est du fait que les sommes sont disjointes dans $Pr(C)$, il 
est facile de montrer le r\'esultat pour le cas o\`u $H$ est isomorphe dans
$Pr(C)$ \`a une somme de pr\'efaisceaux repr\'esentables. 

Venons en au cas g\'en\'eral. Choisissons un \'epimorphisme
$p : X_{0} \longrightarrow H$, avec $X_0$ isomorphe dans $Pr(C)$ \`a une
somme de repr\'esentables.
Soit
$$X_{1}:=X_{0}\times_{H}X_{0} \rightrightarrows X_{0}$$ 
la relation d'\'equivalence sur $X_{0}$ d\'efinie par $p$.
Comme $p$ est un \'epimorphisme dans $Pr(C)$ le morphisme naturel
$$X_{0}/X_{1} \longrightarrow H,$$
du quotient de $X_{0}$ par  la relation $X_{1}$ 
dans $H$, est un isomorphisme dans $Pr(C)$ (cela se voit par exemple
en consid\'erant les fibres de $X_{0}/X_{1} \longrightarrow H$, qui sont
des co-\'egaliseurs des deux projections $Z\times Z \rightrightarrows Z$ pour
un certain ensemble $Z$ non vide). 

Soit $Y_{i}:=f_{!}(X_{i})$ pour $i=0,1$, et
consid\'erons les morphismes naturels
$$Y_{1} \rightrightarrows Y_{0}.$$
Nous d\'eduisons de la proposition \ref{p1} dans le cas o\`u $H$ 
est une somme de repr\'esentables que $Y_{1}$ est encore une relation 
d'\'equivalence sur $Y_{0}$.
De plus, comme $f_!$ est adjoint \`a gauche il commute aux colimites, et on a
$$Y_{0}/Y_{1} \simeq f_!(X_{0}/X_{1}) \simeq f_!(H) \ .$$
Par cons\'equent nous avons 
$$f_!(X_0\times_H X_0) \simeq f_!(X_1) = Y_1 \simeq   
Y_0 \times _{f_{!}(H)}Y_0 \simeq f_!(X_0) \times_{f_!(H)}f_!(X_0) \ .$$ 
Comme le morphisme $p:X_0 \to H$ est un \'epimorphisme, 
les morphismes ${h}_x \to H$ et ${h}_y \to H$ 
se factorisent par $X_0$. Par cons\'equent nous
avons
$${h}_x \times_H {h}_y \simeq {h}_x \times_{X_0}
 ( X_0\times _H X_0) \times_{X_{0}} {h}_y.$$
D'apr\`es  ce qui pr\'ec\`ede cela implique que le morphisme naturel
$$f_!({h}_x \times_H {h}_y) \longrightarrow
f_!({h}_x) \times_{f_!(H)} f_!({h}_y) $$
est un isomorphisme dans $Pr(D)$.  
\hfill $\Box$ \\

Nous supposons maintenant que la cat\'egorie $C$ est munie d'une topologie 
de Grothendieck $\tau$.  La sous-cat\'egorie de $Pr(C)$ form\'ee des faisceaux
pour cette topologie sera not\'ee $Sh(C)$, et le foncteur de faisceautisation 
sera not\'e
$$a : Pr(C) \longrightarrow Sh(C).$$
On rappelle que le foncteur $a$ est exact \`a gauche. 

\begin{prop}\label{p2}
Soient $C$ et $D$ deux cat\'egories poss\'edant des limites finies, et 
$\tau$ (resp. $\rho$) une topologie sur $C$ (resp. sur $D$). Soit
$f : C \longrightarrow D$ un foncteur v\'erifiant les conditions suivantes.
\begin{enumerate}
\item Pour toute famille couvrant $\{U_{i} \longrightarrow X\}$ dans $C$, la
famille $\{f(U_{i}) \longrightarrow f(X)\}$ est couvrante dans $D$.
\item Le foncteur $f$ est exact \`a gauche.
\end{enumerate}
Alors le foncteur 
$$f^{*} : Pr(D) \longrightarrow Pr(C)$$
pr\'eserve les faisceaux. De plus, le foncteur induit
$$f^{*} : Sh(D) \longrightarrow Sh(C)$$
poss\`ede un adjoint \`a gauche
$$f_{!}^{a}:=a\circ f_{!} : \xymatrix{Sh(C) \ar[r]^-{f_{!}} & Pr(D) \ar[r]^-{a} & Sh(D)}$$
qui est exact \`a gauche.
\end{prop}  

\textit{Preuve:} Le fait que $f_{!}^{a}$ d\'efini comme le foncteur compos\'e
$a\circ f_{!}$ soit exact \`a gauche d\'ecoule du fait que $a$ est exact \`a gauche et 
de la proposition \ref{p1}. De plus, si l'on sait que $f^{*}$ pr\'eserve les faisceaux il est formel
de voir que $f_ {!}^{a}$ est adjoint \`a gauche de $f^{*}$. Il nous suffit donc de montrer que
$f^{*}$ pr\'eserve les faisceaux. 

Pour cela, soit $F\in Sh(D)$, $X\in C$ et $\{U_{i} \longrightarrow X\}$ une famille
couvrante. On consid\`ere le morphisme de pr\'efaisceaux
$$U:=\coprod_{i}h_{U_{i}} \longrightarrow h_{X}$$
et $R:=U\times_{h_{X}}U \rightrightarrows U$ la relation d'\'equivalence
d\'efinie sur $U$. Dire que $f^{*}(F)$ poss\`ede la propri\'et\'e de descente
par rapport au recouvrement $\{U_{i} \longrightarrow X\}$ est \'equivalent \`a dire que
le morphisme naturel
$$Hom_{Pr(C)}(U/R,f^{*}(F)) \longrightarrow Hom_{Pr(C)}(h_{X},f^{*}(F))$$
est bijectif. Or, par adjonction ce morphisme est isomorphe au morphisme naturel
$$Hom_{Pr(D)}(f_{!}(U/R),F) \longrightarrow Hom_{Pr(D)}(f_{!}(h_{X}),F).$$
Par hypoth\`ese sur $f$, le quotient $f_{!}(U/R)$ est isomorphe
au quotient $V/R'$, o\`u $V:=\coprod_{i}h_{f(U_{i})}$, et 
$$R':=V\times_{h_{f(X)}}V \rightrightarrows V$$
est la relation d'\'equivalence d\'efinie par la projection $V \longrightarrow h_{f(X)}$.
Ainsi, dire que le morphisme 
$$Hom_{Pr(D)}(f_{!}(U/R),F) \longrightarrow Hom_{Pr(D)}(f_{!}(h_{X}),F)$$
est bijectif est \'equivalent \`a dire que le pr\'efaisceau $F$ poss\`ede la condition
de descente pour la famille couvrante $\{f(U_{i}) \longrightarrow f(X)\}$.
Comme $F$ est un faisceau, cette condition de descente est satisfaite, 
et donc $f^{*}(F)$ poss\`ede la condition de descente pour
$\{U_{i} \longrightarrow X\}$. Cela finit de montrer que $f^{*}(F)$ est un faisceau.
\hfill $\Box$ \\

\subsection{Contextes et champs g\'eom\'etriques}

\begin{df}\label{d1}
Un \emph{contexte g\'eom\'etrique (complet)} est un triplet
$(C,\tau,\mathbb{P})$, o\`u $C$ est une cat\'egorie, $\tau$ est une
topologie de Grothendieck sur $C$ et $\mathbb{P}$ est une classe 
de morphismes dans $C$, v\'erifiant les conditions suivantes. 
\begin{enumerate}
\item La topologie $\tau$ est sous-canonique (nous identifierons alors
$C$ \`a une sous-cat\'egorie pleine de la cat\'egorie $Sh(C)$ des faisceaux sur $C$). 
\item La cat\'egorie $C$ poss\`ede des limites finies et des sommes finies. De plus, les sommes
sont disjointes dans $C$. 
\item La classe de morphismes $\mathbb{P}$ est stable par composition, changements de bases, et contient
les isomorphismes. 
\item Pour toute famille couvrante $\{X_{i}\longrightarrow X\}$ dans $C$, 
il existe des morphismes $Y_{i} \longrightarrow X$ dans $\mathbb{P}$ et des diagrammes commutatifs
$$\xymatrix{
 & X_{i} \ar[d] \\
Y_{i} \ar[ru]\ar[r] & X}$$ 
tels que nous ayons un \'epimorphisme de faisceaux 
$$\coprod_{i}h_{Y_{i}} \longrightarrow h_{X}$$
\item Soit $f : Y \longrightarrow X$ un morphisme dans $C$ tel qu'il existe une famille couvrante $\{Y_{i} \longrightarrow Y\}$ dans $C$, telle que
tous les morphismes $Y_{i} \longrightarrow Y$ et $Y_{i} \longrightarrow X$ soient dans $\mathbb{P}$. Alors le morphisme
$f$ est dans $\mathbb{P}$.
\item La cat\'egorie $C$ poss\`ede des sommes finies, et pour toute paire d'objets $X$ et $Y$ dans $C$ le morphisme naturel
$$X \longrightarrow X\coprod Y$$
est dans $\mathbb{P}$. De plus, la famille de morphismes
$$\{X\longrightarrow X\coprod Y, Y \longrightarrow X\coprod Y\}$$
est couvrante.
\item Soit $h : C \longrightarrow Sh(C)$, le plongement de Yoneda de $C$ dans la cat\'egorie des faisceaux sur $C$ pour 
la topologie $\tau$. Soit $X$ un objet de $C$ tel qu'il existe un isomorphisme dans $Sh(C)$
$$h_{X}\simeq F_{1}\coprod F_{2}.$$
Alors il existe des objets $X_{1}$ et $X_{2}$ dans $C$ tels 
$$F_{1}\simeq h_{X_{1}} \qquad F_{2}\simeq h_{X_{2}}.$$
\item Soient $X\in C$ et $F \longrightarrow X$ un morphisme de faisceaux sur $C$. On suppose qu'il existe
une famille couvrante $\{X_{i} \longrightarrow X\}$ dans $C$ telle que chacun des faisceaux
$F\times_{X}X_{i}$ soit repr\'esentable par un objet de $C$. Alors $F$ est repr\'esentable par un objet de $C$.
\end{enumerate}
\end{df}

\begin{rmk}\label{r1}
\emph{Nous ne mentionerons pas le terme} complet \emph{de la d\'efinition \ref{d1}. Il fait r\'ef\'erence aux
propri\'et\'es $(2)$ et $(8)$, alors qu'il existe une notion plus g\'en\'erale de contextes
g\'eom\'etrique pour la quelle ces conditions ne sont pas satisfaites (voir \cite{to-dea}). Tous les contextes
g\'eom\'etriques que nous utiliserons dans ce travail seront complets par convention}.
\end{rmk}

Fixons un contexte g\'eom\'etrique $(C,\tau,\mathbb{P})$, et consid\'erons la
cat\'egorie $Sh(C)$ des faisceaux sur le site $(C,\tau)$. 
Nous allons d\'efinir une sous-cat\'egorie pleine de $Sh(C)$ form\'ee
des \emph{espaces g\'eom\'etriques} obtenus par recollement 
d'objets de $C$. La d\'efinition se fait par r\'ecurrence et suit les m\^emes
\'etapes que celle de \cite{hagII}. Pour cela rappelons que l'on dispose d'un plongement de Yoneda
$$h : C \longrightarrow Pr(C)$$
qui se factorise en un foncteur pleinement fid\`ele (car $\tau$ est suppos\'ee sous-canonique)
$$h : C \longrightarrow Sh(C).$$
Nous identifierons la cat\'egorie $C$ \`a son image essentielle par le foncteur $h$, et 
omettrons ainsi de mentioner $h$.

\begin{enumerate}
\item Nous commen\c{c}ons par d\'ecr\'eter qu'un faisceau
$F \in Sh(C)$ est \emph{(-1)-g\'eom\'etrique} s'il est 
repr\'esentable (i.e. isomorphe dans $Sh(C)$ \`a un objet
de la forme $h_{X}$ pour $X\in C$).

\item Un morphisme $f : F \longrightarrow G$ de faisceaux est 
\emph{(-1)-repr\'esentable} si pour tout 
$X\in C$, et tout morphisme $X \longrightarrow G$
le faisceau $F\times_{G}X$ est (-1)-g\'eom\'etrique.

\item Un morphisme (-1)-repr\'esentable $f : F \longrightarrow G$ 
est \emph{dans $\mathbb{P}$} si pour tout $X\in C$ et tout morphisme
$X \longrightarrow G$ le morphisme induit entre faisceaux repr\'esentables
$$F\times_{G}X \longrightarrow X$$
est dans $\mathbb{P}$. 

\item Soit $n\geq 0$ et supposons que les trois notions de faisceaux $m$-g\'eom\'etriques,
de morphismes $m$-repr\'esentables et de morphismes $m$-repr\'esentables
et dans $\mathbb{P}$ soient d\'efinies pour tout $m<n$. 
\begin{enumerate}
\item Un faisceau
$F$ est \emph{$n$-g\'eom\'etrique} s'il existe une famille d'objets
$\{U_{i}\}$ dans $C$ et un morphisme
$$U:=\coprod_{i}U_{i} \longrightarrow F$$ 
satisfaisant les deux conditions suivantes
\begin{enumerate}
\item Le morphisme $U \longrightarrow F$ est un \'epimorphisme 
de faisceaux. 
\item Chacun des morphismes $U_{i} \longrightarrow F$
est $(n-1)$-repr\'esentable et dans $\mathbb{P}$. 
\end{enumerate}
La donn\'ee des $U_{i}$ et du morphisme
$$\coprod_{i}U_{i} \longrightarrow F$$
sera appel\'ee un \emph{n-atlas pour $F$}.
\item  Un morphisme $f : F \longrightarrow G$ est \emph{$n$-repr\'esentable}
si pour tout $X\in C$ et tout morphisme $X\longrightarrow G$ le faisceau
$F\times_{G}X$ est $n$-g\'eom\'etrique.
\item Un morphisme $n$-repr\'esentable $f : F \longrightarrow G$ 
est \emph{dans $\mathbb{P}$} si pour tout 
$X\in C$ et tout morphisme $X \longrightarrow G$ il existe un
$n$-atlas $\{U_{i}\}$ de $F\times_{G}X$ tel que tous les morphismes 
induits 
$$U_{i} \longrightarrow X$$
soient dans $\mathbb{P}$. 
\end{enumerate}
\end{enumerate}

Nous posons alors la d\'efinition suivante. 

\begin{df}\label{d2}
\begin{enumerate}
\item
Un faisceau $F\in Sh(C)$ est \emph{g\'eom\'etrique} s'il est $n$-g\'eom\'etrique pour un 
certain $n$. 
\item Un morphisme de faisceaux $f : F \longrightarrow G$ est \emph{dans $\mathbb{P}$} s'il est 
$n$-repr\'esentable et dans $\mathbb{P}$ pour un certain $n$.
\end{enumerate}
\end{df}

On remarque que le proc\'ed\'e inductif de la d\'efinition pr\'ec\'edente s'arr\`ete
en r\'ealit\'e \`a $n=1$. 

\begin{lem}\label{l1}
Si $n>1$, alors tout faisceau $n$-g\'eom\'etrique est $(n-1)$-g\'eom\'etrique.
\end{lem}

\textit{Preuve:} Soit $F$ un faisceau $n$-g\'eom\'etrique et 
$$U:=\coprod_{i}U_{i} \longrightarrow F$$
un $n$-atlas. Nous allons montrer que les morphismes 
$U_{i} \longrightarrow F$ sont en fait $0$-g\'eom\'etrique. Pour cela, soit $X\in C$ et 
$X\longrightarrow F$ un morphisme. Le faisceau $U_{i}\times_{F}X$ est un sous-faisceau
de $U_{i}\times X$. Soit 
$$\coprod_{j}V_{i,j} \longrightarrow U_{i}\times_{F}X$$
un $(n-1)$-atlas. Alors, pour tout $Y\in C$ et tout morphisme $Y\longrightarrow U_{i}\times_{F}X$, on a
$$V_{i,j}\times_{U_{i}\times_{F}X}Y\simeq V_{i,j}\times_{U_{i}\times X}Y.$$
Comme $C$ a des limites finies on voit que $V_{i,j}\times_{U_{i}\times_{F}X}Y$ est 
repr\'esentable pour tout $j$ et tout $i$. Cela implique que pour tout $i$ et tout $j$, le morphisme
$$V_{i,j} \longrightarrow U_{i}\times_{F}X$$
est $(-1)$-g\'eom\'etrique, et donc que 
$U_{i}\times_{F}X$ est $0$-g\'eom\'etrique.
\hfill $\Box$ \\

\begin{rmk}\label{r2}
\emph{En termes g\'eom\'etriques, les faisceaux $(-1)$-g\'eom\'etriques correspondent
aux objets affines, les faisceaux $0$-g\'eom\'etriques aux espaces ayant une
diagonale affine, et les faisceaux $1$-g\'eom\'etriques aux espaces sans
aucune condition de s\'eparation.}
\end{rmk}

Comme dans \cite{hagII} on v\'erifie que les faisceaux g\'eom\'etriques sont stables par 
limites finies, sommes disjointes infinies et quotients par des
relations d'\'equivalences qui sont dans $\mathbb{P}$. On renvoie \`a \cite{hagII} pour plus de d\'etails, que nous
ne reproduirons pas ici. Nous utiliserons en particulier le lemme suivant. \\

\begin{lem}\label{l2}
Soit $X\in C$ et $F \longrightarrow X$ un morphisme dans $Sh(C)$. 
\begin{enumerate}
\item
Supposons qu'il
existe une famille couvrante $\{X_{i} \longrightarrow X\}$ dans $C$ telle que pour tout
$i$ le faisceau $F\times_{X}X_{i}$ soit g\'eom\'etrique. Alors le faisceau
$F$ est g\'eom\'etrique. 
\item Supposons qu'il
existe une famille couvrante $\{X_{i} \longrightarrow X\}$ dans $C$ telle que pour tout
$i$ le morphisme $F\times_{X}X_{i} \longrightarrow X_{i}$ soit 
g\'eom\'etrique et dans $\mathbb{P}$. Alors le faisceau
$F$ est g\'eom\'etrique et $F \longrightarrow X$ est dans $\mathbb{P}$.
\end{enumerate}
\end{lem}

\textit{Preuve:} D'apr\`es le point $(4)$ de la d\'efinition \ref{d1} on peut supposer que 
pour tout $i$ le morphisme $X_{i} \longrightarrow X$ est dans $\mathbb{P}$. Si pour tout
$i$ on choisit un $1$-atlas
$\{U_{i,j} \longrightarrow F\times_{X}X_{i}\}$, alors on v\'erifie facilement que la famille totale
$\{U_{i,j} \longrightarrow F\}$ est encore un $1$-atlas. Ceci implique que $F$ est g\'eom\'etrique.
Si de plus chaque $F\times_{X}X_{i} \longrightarrow X_{i}$ est dans $\mathbb{P}$, alors
on peut choisir $U_{i,j}$ de sorte \`a ce que $U_{i,j} \longrightarrow X_{i}$
soit dans $\mathbb{P}$. Ceci implique
que $F \longrightarrow X$ est dans $\mathbb{P}$.
\hfill $\Box$ \\

Nous terminons cette section par la notion de changement de contextes g\'eom\'etriques. Pour cela, soient
$(C,\tau,\mathbb{P})$ et $(D,\rho,\mathbb{Q})$ deux contextes g\'eom\'etriques. Soit 
$f : C \longrightarrow D$ un foncteur.

\begin{prop}\label{p3}
On suppose que le foncteur $f$ v\'erifie les deux conditions suivantes.
\begin{enumerate}
\item Le foncteur $f$ commute aux limites finies et est continu pour les topologies $\tau$ et $\rho$
(i.e. $f^{*}$ pr\'eserve les faisceaux). 
\item L'image par $f$ d'un morphisme de $\mathbb{P}$ est un morphisme 
de $\mathbb{Q}$. 
\end{enumerate}
Alors le foncteur induit de la proposition \ref{p2}
$$f^{a}_{!} : Sh(C) \longrightarrow Sh(D)$$
pr\'eserve les faisceaux g\'eom\'etriques et envoie les morphismes
de $\mathbb{P}$ dans des morphismes de $\mathbb{Q}$. 
\end{prop}

\textit{Preuve:} Commen\c{c}ons par montrer que 
$f^{a}_{!}$ pr\'eserve les morphismes $n$-repr\'esentables et envoie
les morphismes $n$-repr\'esentables et dans $\mathbb{P}$ dans 
des morphismes $n$-repr\'esentables et dans $\mathbb{Q}$. Nous proc\`ederons par r\'ecurrence sur $n$. 
Pour $n=-1$ c'est la formule 
$$f^{a}_{!}(h_{X})\simeq h_{f(X)}$$
et les conditions sur le foncteur $f$.
Supposons que tout cela soit vrai pour $m<n$, et soit 
$g : F \longrightarrow G$ un morphisme $n$-repr\'esentable dans $Sh(C)$. Soit $Y\in D$ et 
$Y\longrightarrow f^{a}_{!}(G)$ un morphisme dans $Sh(D)$. 

\begin{lem}\label{l3}
Il existe un famille d'objets $\{X_{i}\}_{i\in I}$ dans $C$, une 
famille couvrante $\{Y_{j} \longrightarrow Y\}_{j\in J}$ dans $D$,
une application $u : J \rightarrow I$, un morphisme
$$\coprod_{i}X_{i} \longrightarrow G$$
et pour tout $j$ un diagramme commutatif dans $Sh(D)$
$$\xymatrix{
Y_{j}  \ar[r] \ar[d] & f^{a}_{!}(X_{u(j)})\simeq f(X_{u(j)}) \ar[d] \\
Y \ar[r]  & f^{a}_{!}(G).
}$$
\end{lem}

\textit{Preuve du lemme \ref{l3}:} Le foncteur 
$f^{a}_{!}$ pr\'eserve les \'epimorphismes de faisceaux, car il commute 
aux limites finies et aux colimites arbitraires.
Soit donc $\coprod_{i}X_{i} \longrightarrow G$ un \'epimorphisme dans $Sh(C)$
avec $X_{i}$ des objets de $C$. Alors le morphisme induit
$$\coprod_{i}f(X_{i}) \longrightarrow f^{a}_{!}(G)$$
est un \'epimorphisme dans $Sh(D)$, ce qui implique le lemme par d\'efinition des
\'epimorphismes. 
\hfill $\Box$ \\

Revenons \`a la preuve de la proposition \ref{p3}. Il nous faut montrer que 
$f^{a}_{!}(F)\times_{f^{a}_{!}(G)}Y$ est un faisceau 
$n$-g\'eom\'etrique. Le lemme pr\'ec\'edent et le lemme \ref{l2}  permettent de supposer que 
le morphisme $Y \longrightarrow f^{a}_{!}(G)$ se factorise par un 
morphisme $f^{a}_{!}(X) \longrightarrow f^{a}_{!}(G)$ pour un $X\in C$ et un morphisme
$Y \longrightarrow f^{a}_{!}(X)$. Mais dans ce cas, on a
$$f^{a}_{!}(F)\times_{f^{a}_{!}(G)}Y\simeq
f^{a}_{!}(F)\times_{f^{a}_{!}(G)}f^{a}_{!}(X)\times_{f^{a}_{!}(X)}Y\simeq
f^{a}_{!}(F\times_{G}^{h}X)\times_{f(X)}Y.$$
On est ainsi ramen\'es \`a d\'emontrer que $f^{a}_{!}$ pr\'eserve
les faisceaux $n$-g\'eom\'etriques, mais cela se d\'eduit imm\'ediatement de l'hypoth\`ese de r\'ecurrence.

Il nous reste \`a voir que si de plus le morphisme $g$ est dans $\mathbb{P}$ alors
$f^{a}_{!}(g)$ est dans $\mathbb{Q}$. Le m\^eme argument montre alors
que l'on peut supposer que $G=X$ est un objet de $C$. Dans ce cas cela se d\'eduit encore
des hypoth\`eses de r\'ecurrence.
\hfill $\Box$ \\

\subsection{Analytification des espaces alg\'ebriques}

Nous allons maintenant sp\'ecifier deux contextes g\'eom\'etriques, le contexte
alg\'ebrique et le contexte analytique (tous deux au-dessus des nombres complexes). \\

Nous commencerons par le contexte alg\'ebrique. Pour cela nous
d\'efinissons un contexte g\'eom\'etrique $(C,\tau,\mathbb{P})$, pour lequel
$C:=Aff$ est la cat\'egorie des $\mathbb{C}$-sch\'emas affines
et de type fini, $\tau:=et$ est la topologie \'etale, et $\mathbb{P}$ est la classe
des morphismes \'etales. Ces notions d\'efinissent un contexte g\'eom\'etrique
$(Aff,et,et)$
au sens de la d\'efinition \ref{d1}. 

\begin{df}\label{d3}
\begin{enumerate}
\item
Un \emph{espace alg\'ebrique} 
est un faisceau g\'eom\'etrique pour
le contexte $(Aff,et,et)$ d\'efini ci-dessus.
\item La sous-cat\'egorie pleine de $Sh(Aff)$ form\'ee des espaces
alg\'ebriques sera not\'ee
$Esp^{alg}$.
\end{enumerate}
\end{df}

\begin{rmk}\label{r3}
\emph{Tels que nous les avons d\'efinis ci-dessus les espaces alg\'ebriques seront toujours
localement de type fini sur $\mathbb{C}$. Cela provient du choix du contexte
$(Aff,et,et)$ pour lequel nous nous restreignons aux sch\'emas affines
de type fini. Il existe une notion plus g\'en\'erale d'espaces alg\'ebriques non-n\'ecessairement
localement de type fini en consid\'erant tous les sch\'emas affines (voir \cite{hagII}). 
Cette notion plus g\'en\'erale ne sera
pas utilis\'ee dans ce travail.}
\end{rmk}

Passons maintenant \`a la description du contexte analytique $(C,\tau,\mathbb{P})$.
Pour cela nous posons $C:=Ste$, la cat\'egorie des espaces analytiques de Stein\footnote{Dans ce travail, 
un espace analytique $X$ est de Stein s'il est d\'enombrable \`a l'infini et si
$H^{i}(X,\mathcal{F})=0$ pour tout $i>0$ et tout faisceau coh\'erent $\mathcal{F}$.}, 
$\tau:=top$ la topologie usuelle sur les espaces analytiques complexes, et 
pour $\mathbb{P}$ nous prenons la classe des morphismes \'etales entre espaces
analytiques (i.e. des isomorphismes biholomorphes locaux).
Ces notions d\'efinissent un contexte g\'eom\'etrique $(Ste,top,et)$ au sens de la d\'efinition \ref{d1}.

\begin{df}\label{d4}
\begin{enumerate}
\item
Un \emph{espace analytique} 
est un faisceau g\'eom\'etrique pour
le contexte $(Ste,top,et)$ d\'efini ci-dessus.
\item La sous-cat\'egorie pleine de $Sh(Ste)$ form\'ee des espaces
analytiques sera not\'ee
$Esp^{an}$.
\end{enumerate}
\end{df}

\begin{rmk}\label{r4}
\emph{Les d\'efinitions d'espaces alg\'ebriques et d'espaces analytiques ci-dessus sont plus g\'en\'erales que celles
que l'on rencontre usuellement dans la litt\'erature. Par exemple, nos espaces alg\'ebriques n'ont pas des diagonales
quasi-compactes alors que cette hypoth\`ese est souvent incluse dans la d\'efinition d'espaces
alg\'ebriques (cf. \cite{Kn}). De m\^eme, nos espaces analytiques sont plus g\'en\'eraux que les espaces
analytiques usuels. Par exemple, si $\Gamma$ est un groupe discret qui op\`ere 
sans points fixes, mais non proprement, sur un espace analytique $X$, alors le faisceau quotient
$X/\Gamma$ est un espace analytique au sens de la d\'efinition pr\'ec\'edente. Cependant, ce quotient peut ne pas
exister en tant qu'espace analytique au sens usuel (e.g. au sens de \cite{bs,gr}) . De plus, 
m\^eme lorsqu'un quotient $X//\Gamma$ existe 
en tant qu'espace analytique au sens usuel, en g\'en\'eral $X/\Gamma$ et $X//\Gamma$ ne sont pas isomorphes.
Le lecteur pourra garder en t\^ete l'exemple suivant: $X=\mathbb{C}^{\times}$, et $\Gamma=\mathbb{Z}$ qui 
op\`ere par $z\mapsto q.z$ avec $|q|=1$ et qui n'est pas une racine de l'unit\'e. Dans ce cas
$X//\Gamma$ est r\'eduit \`a un point mais $X/\Gamma$ ne l'est pas.}
\end{rmk}

Consid\'erons maintenant le foncteur d'analytification (voir \cite[Exp. XII]{sga1})
$$\begin{array}{ccc}
Aff & \longrightarrow & Ste \\
X & \mapsto & X^{an}.
\end{array}$$
Ce foncteur poss\`ede toutes les propri\'et\'es de la proposition \ref{p3}, et ainsi 
donne lieu \`a un foncteur au niveau des faisceaux
$$(-)^{an} : Sh(Aff) \longrightarrow Sh(Ste).$$
Ce foncteur commute aux limites finies, envoie
les espaces alg\'ebriques sur les espaces analytiques et
pr\'eserve les morphismes \'etales. Nous disposons donc d'un foncteur induit
$$(-)^{an} : Esp^{alg} \longrightarrow Esp^{an}.$$

Pour terminer cette section nous allons donner une description explicite 
de l'analytifi\'e d'un faisceaux $F\in Sh(Aff)$, tout au moins localement 
en chaque points. Pour cela, si $x\in X$ est un point dans un
espace de Stein et si $F\in Sh(Ste)$ est un faisceau nous noterons
$$F(X_{x}):=Colim_{x\in U\subset X}F(U),$$
o\`u la colimite est prise sur l'ensemble filtrant des ouverts de Stein de $X$ contenant $x$. 
En d'autres termes, $F(X_{x})$ est la fibre du faisceau $F$ restreint au petit
site des ouverts de Stein de $X$. La notation $X_{x}$ d\'esigne 
le germe d'espaces analytiques de $X$ en $x$. Pour tout $Y\in Ste$ nous poserons aussi
$$Hom(X_{x},Y):=Colim_{x\in U\subset X}Hom(U,Y),$$
les germes de morphismes en $x$ de $X$ vers $Y$.

Soit $\mathcal{O}_{x}$ la $\mathbb{C}$-alg\`ebre des germes de fonctions
holomorphes sur $X$ au point $x$. Pour toute $\mathbb{C}$-alg\`ebre
commutative de type fini $B$, il existe une application
$$Hom(X_{x},(Spec\, B)^{an}) \longrightarrow Hom_{\mathbb{C}-Alg}(B,\mathcal{O}_{x}),$$
qui \`a $u : X_{x} \longrightarrow (Spec\, B)^{an}$ associe le morphisme
$$u^{*} : \mathcal{O}((Spec\, B)^{an}) \longrightarrow \mathcal{O}(X_{x})=\mathcal{O}_{x}$$
compos\'e \`a droite avec le morphisme naturel
$$B \longrightarrow \mathcal{O}((Spec\, B)^{an})$$
(qui \`a un \'el\'ement de $B$ associe la fonction holomorphe sur $(Spec\, B)^{an}$
correspondante).

\begin{lem}\label{l4}
L'application ci-dessus
$$Hom(X_{x},(Spec\, B)^{an})\longrightarrow Hom_{\mathbb{C}-Alg}(B,\mathcal{O}_{x}),$$
est bijective.
\end{lem}

\textit{Preuve du lemme:} Il est facile de voir que l'ensemble des $\mathbb{C}$-alg\`ebres $B$
pour lesquelles la conclusion du lemme est v\'erifi\'ee est stable
par colimites finies. Il nous suffit donc de v\'erifier le lemme pour $B=\mathbb{C}[X]$, 
qui est alors \'evident. 
\hfill $\Box$ \\

On consid\`ere maintenant l'ensemble filtrant des sous-$\mathbb{C}$-alg\`ebres
de type fini $A\subset \mathcal{O}_{x}$.
Notons $f : Aff \longrightarrow Ste$ le foncteur d'analytification, et 
$f_{!} : Pr(Aff) \longrightarrow Pr(Ste)$ le foncteur induit sur les pr\'efaisceaux. 
Comme ci-dessus, on pose pour tout $F\in Pr(Ste)$
$$F(X_{x}):=Colim_{x\in U\subset X}F(U).$$

Le lemme pr\'ec\'edent implique que pour
toute sous-$\mathbb{C}$-alg\`ebre de type finie $A\subset \mathcal{O}_{x}$
on dispose d'un germe de morphisme naturel
$$X_{x} \longrightarrow (Spec\, A)^{an}.$$
Pour $F\in Pr(Aff)$, on dispose donc
d'un morphisme naturel 
$$f_{!}(F)((Spec\, A)^{an}) \longrightarrow f_{!}(F)(X_{x}).$$
En composant avec le morphisme d'analytification
$$F(Spec\, A)=Hom(Spec\, A,F) \longrightarrow f_{!}(F)((Spec\, A)^{an})=Hom(f_{!}(Spec\, A),f_{!}(F))$$
on en d\'eduit un morphisme 
$$F(Spec\, A) \longrightarrow f_{!}(F)(X_{x})$$
fonctoriel en $A\subset \mathcal{O}_{x}$. Cela d\'efinit un morphisme
$$\phi : Colim_{A\subset \mathcal{O}_{x}}F(Spec\, A) \longrightarrow f_{!}(F)(X_{x}).$$

\begin{prop}\label{p4}
Avec les notations ci-dessus le morphisme
$$\phi : Colim_{A\subset \mathcal{O}_{x}}F(Spec\, A) \longrightarrow f_{!}(F)(X_{x})$$
est bijectif.
\end{prop}

\textit{Preuve:}
Nous dirons qu'un pr\'efaisceau $F\in Pr(Aff)$ \`a la propri\'et\'e $P$, si 
le morphisme de la proposition
$$\phi : Colim_{A\subset \mathcal{O}_{x}}F(Spec\, A) \longrightarrow f_{!}(F)(X_{x})$$
est bijectif (le couple $(X,x)$ \'etant fix\'e). Il est facile de voir que la sous-cat\'egorie pleine
de $Pr(Aff)$ form\'ee des pr\'efaisceaux ayant la propri\'et\'e $P$
est stable par colimites et limites finies. Comme tout pr\'efaisceau est 
colimite de repr\'esentables, on voit qu'il nous suffit de montrer que 
les repr\'esentables poss\`edent la propri\'et\'e $P$. Soit donc $F=Spec\, B$, avec 
$B$ une $\mathbb{C}$-alg\`ebre de type finie. Le pr\'efaisceau $F$ s'\'ecrit alors
comme un produit fibr\'e de la forme
$$\xymatrix{
F \ar[r] \ar[d] & (Spec\, \mathbb{C}[X])^{n} \ar[d] \\
Spec\, \mathbb{C} \ar[r] & (Spec\, \mathbb{C}[Y])^{m}.}$$
Comme l'ensemble des pr\'efaisceaux poss\`edant la propri\'et\'e $P$ est stable
par limites finies on est ramen\'es au cas o\`u $F=\mathbb{A}^{1}=Spec\, \mathbb{C}[X]$
pour lequel la proposition est \'evidente.
\hfill $\Box$ \\

Supposons maintenant que $F\in Sh(Aff)$ soit un faisceau.
Le morphisme de faisceautisation
$$f_{!}(F) \longrightarrow a(f_{!}(F))=F^{an}$$
induit alors une bijection
$f_{!}(X_{x}) \simeq F^{an}(X_{x})$. Le morphisme $\phi$ de la proposition \ref{p4} peut donc aussi 
\^etre consid\'er\'e comme une bijection
$$\phi : Colim_{A\subset \mathcal{O}_{x}}F(Spec\, A) \simeq F^{an}(X_{x}).$$
C'est sous cette derni\`ere forme que nous utiliserons la proposition \ref{p4}.

\section{Cat\'egories d\'eriv\'ees des espaces analytiques}

Pour un espace analytique $X$ (voir nos conventions \`a la section \S 2)
pour la notion d'espace analytique que nous utilisons), nous disposons de sa cat\'egorie
des $\mathcal{O}_{X}$-modules $Mod(\mathcal{O}_{X})$. C'est une cat\'egorie
de Grothendieck, et nous noterons $D(X)$ sa cat\'egorie d\'eriv\'ee non-born\'ee
(voir l'appendice C).
La sous-cat\'egorie pleine
de $D(X)$ form\'ee des complexes \`a cohomologie coh\'erente
sera not\'ee $D_{coh}(X)$. 

\begin{df}\label{d5}
Un objet $E\in D(X)$ est \emph{parfait} (nous dirons aussi 
\emph{complexe parfait}) s'il est localement sur $X$ 
isomorphe \`a un complexe born\'e de $\mathcal{O}_{X}$-modules localement
libres de rangs finis. 

La sous-cat\'egorie pleine des objets parfaits dans $D(X)$ sera not\'ee
$D_{parf}(X)$.
\end{df}

Nous avons bien entendu des inclusions
$$D_{parf}(X) \subset D_{coh}(X) \subset D(X).$$

La cat\'egorie $D(X)$ est munie d'une structure tensorielle 
$-\otimes^{\mathbb{L}}-$ (voir l'appendice C). La sous-cat\'egorie
$D_{parf}(X)$ est stable par $-\otimes^{\mathbb{L}}-$ et h\'erite donc
elle aussi d'une structure tensorielle. \\

Soit maintenant $f : X \longrightarrow Y$ un morphisme d'espaces analytiques. 
Il induit une adjonction 
$$f^{*} : Mod(\mathcal{O}_{Y}) \rightleftarrows Mod(\mathcal{O}_{X}) : f_{*}.$$
Comme il est expliqu\'e dans l'appendice C
cette adjonction se d\'erive alors en une adjonction au niveau des cat\'egories
d\'eriv\'ees
$$\mathbb{L}f^{*} : D(Y) \rightleftarrows D(X) : \mathbb{R}f_{*}.$$
Le foncteur $\mathbb{L}f^{*}$ est naturellement compatible avec 
les structures tensorielles.

\begin{prop}\label{p5}
Soit $f : X \longrightarrow Y$ un morphisme propre et plat 
d'espaces analytiques. Alors $\mathbb{L}f^{*}$ et 
$\mathbb{R}f_{*}$ pr\'eservent les objets parfaits et induisent une
adjonction
$$\mathbb{L}f^{*} : D_{parf}(Y) \rightleftarrows D_{parf}(X) : \mathbb{R}f_{*}.$$
\end{prop}

\textit{Preuve:} Que $\mathbb{L}f^{*}$ pr\'eserve les parfaits est clair, et ne
demande d'ailleurs pas que $f$ soit propre ou plat. Il nous faut donc
montrer que $\mathbb{R}f_{*}$ pr\'eserve les complexes parfaits.
Nous allons d\'eduire ceci des \'enonc\'es de
finitude et de changement de bases pour les faisceaux analytiques coh\'erents.

\begin{lem}\label{l5}
Pour tout complexe parfait $E$ sur $X$ le complexe
$\mathbb{R}f_{*}(E)$ est \`a cohomologie coh\'erente et localement born\'ee sur
$Y$. 
\end{lem}

\textit{Preuve du lemme \ref{l5}:} Pour $E$ un faisceau coh\'erent le lemme
se d\'eduit de \cite[Thm. 4.1]{bs}(prendre $\mathcal{M}=\mathcal{O}_{Y}$). 
Par r\'ecurrence sur le nombre de faisceaux de cohomologie non-nuls on d\'eduit alors
le lemme pour tout complexe $E$ \`a cohomologie 
coh\'erente et born\'ee. Par propret\'e on en d\'eduit le lemme pour
tout complexe $E$ \`a cohomologie coh\'erente et localement born\'ee.
En particulier il est vrai pour les complexes parfaits. 
\hfill $\Box$ \\

\begin{lem}\label{l6}
Pour tout complexe parfait $E$ sur $X$, et tout 
faisceau coh\'erent $F$ sur $Y$, le morphisme naturel
$$\mathbb{R}f_{*}(E)\otimes^{\mathbb{L}}F \longrightarrow 
\mathbb{R}f_{*}(E\otimes^{\mathbb{L}} \mathbb{L}f^{*}(F))$$ 
est un isomorphisme sur $Y$.
\end{lem}

\textit{Preuve du lemme \ref{l6}:} Comme pour le lemme \ref{l5} on se ram\`ene \`a d\'emontrer
l'assertion du lemme pour  un faisceau coh\'erent $E$ sur $X$. De plus, 
l'assertion \'etant locale sur $Y$, on peut supposer que 
le faisceau $F$ poss\`ede une r\'esolution par des $\mathcal{O}_{Y}$-libres
de rangs finis
$$\xymatrix{\dots & 
F_{n} \ar[r] & F_{n-1} \ar[r] & \dots & F_{1} \ar[r] & F_{0} \ar[r] & F}$$
(on peut par exemple remplacer $Y$ par un recouvrement 
par des ouverts de Stein relativement compacts dans des voisinages ouverts de Stein).
On peut aussi supposer d'apr\`es le lemme \ref{l5} que $\mathbb{R}f_{*}(E)$
est \`a cohomologie born\'ee sur $Y$. 

Soit $p\in \mathbb{Z}$. 
Nous affirmons qu'il existe un entier $m$, tel que le morphisme naturel
$$F_{*}[\leq m] \longrightarrow F$$
induise un isomorphisme sur les faisceaux de cohomologie
$$\underline{H}^{p}(\mathbb{R}f_{*}(E)\otimes F_{*}[\leq m]) \longrightarrow 
\underline{H}^{p}(\mathbb{R}f_{*}(E)\otimes F_{*})$$
$$\underline{H}^{p}(\mathbb{R}f_{*}(E\otimes f^{*}(F_{*}[\leq m]))) \longrightarrow 
\underline{H}^{p}(\mathbb{R}f_{*}(E\otimes f^{*}(F_{*}))),$$
o\`u $F_{*}[\leq m]$ est la r\'esolution tronqu\'ee d\'efinie par
$$(F_{*}[\leq m])_{n}=F_{n} \; si \; n\leq m\; et \; 0 \; sinon.$$
Ces deux assertions se d\'eduisent ais\'ement du fait que
$\mathbb{R}f_{*}(E)$ soit un complexe \`a cohomologie born\'ee sur $Y$. 
Le diagramme commutatif dans $D(Y)$ 
$$\xymatrix{
\mathbb{R}f_{*}(E\otimes f^{*}(F_{*}[\leq m])) \ar[r] \ar[d] & \mathbb{R}f_{*}(E\otimes f^{*}(F_{*})) \ar[d] \\
\mathbb{R}f_{*}(E)\otimes F_{*}[\leq m] \ar[r] & \mathbb{R}f_{*}(E)\otimes F_{*}}$$
montre alors qu'il suffit de d\'emontrer que pour tout $m$ le morphisme naturel
$$\mathbb{R}f_{*}(E\otimes f^{*}(F_{*}[\leq m])) \longrightarrow \mathbb{R}f_{*}(E)\otimes F_{*}[\leq m] $$
est un isomorphisme dans $D(Y)$. Mais cela se ram\`ene facilement par r\'ecurrence sur $m$ 
au cas o\`u $m=0$ qui est \'evident.
\hfill $\Box$ \\

Nous pouvons maintenant d\'emontrer la proposition \ref{p5}. 
Comme l'assertion est locale sur $Y$ on peut utiliser le lemme \ref{l5} et 
supposer que $\mathbb{R}f_{*}(E)$ est \`a cohomologie born\'ee sur $Y$. On applique
alors le lemme \ref{l6} \`a $F=k(y)$, le faisceau gratte-ciel en un point
$y\in Y$. Comme $f$ est plat on a $\mathbb{L}f^{*}(k(y))\simeq \mathcal{O}_{X_{y}}$, 
o\`u $X_{y}$ est la fibre de $f$ en $y$.

Cela nous dit que le morphisme naturel
$$\mathbb{R}f_{*}(E)\otimes^{\mathbb{L}}k(y) \longrightarrow 
\mathbb{R}f_{*}(E\otimes^{\mathbb{L}}\mathcal{O}_{X_{y}})$$
est un quasi-isomorphisme. Or, le second membre est quasi-isomorphe
au complexe d'hyper-cohomologie de $X_{y}$ \`a coefficient dans
le complexe parfait $E\otimes^{\mathbb{L}}\mathcal{O}_{X_{y}}$. Ceci implique
donc que pour tout $y\in Y$, le complexe
$\mathbb{R}f_{*}(E)\otimes^{\mathbb{L}}k(y)$ est cohomologiquement 
born\'e. 

\begin{lem}\label{l7}
Soit $F$ un complexe sur $Y$ \`a cohomologie 
coh\'erente et born\'ee. Soit $y\in Y$ tel que le complexe
$F\otimes^{\mathbb{L}}k(y)$ soit cohomologiquement born\'e. Alors il 
existe un voisinage $Y'$ de $y$ tel que la restriction de $F$ sur $Y'$
soit un complexe parfait. 
\end{lem}

\textit{Preuve du lemme \ref{l7}:} Quitte \`a restreindre $F$ \`a un voisinage
de $y$ on peut supposer que $F$ est quasi-isomorphe \`a un complexe $F'_{*}$
de $\mathcal{O}_{Y}$-modules libres de rangs finis born\'e \`a droite. 
Supposons alors que $H^{i}(F\otimes^{\mathbb{L}}k(y))=0$
pour tout $i<n$, pour $n\leq 0$. Soit 
$\tau_{\geq n}F'_{*}$ le tronqu\'e de $F'_{*}$ 
(par d\'efinition la cohomologie de $\tau_{\geq n}F'_{*}$ coincide
avec celle de $F'_{*}$ en degr\'es sup\'erieurs \`a $n$ et est nulle 
ailleurs). On a donc
$$(\tau_{\geq n}F'_{*})_{n}:=Coker(F'_{n-1} \longrightarrow F'_{n}).$$
Par hypoth\`ese on trouve que 
$$Tor_{i}^{\mathcal{O}_{Y}}((\tau_{\geq n}F'_{*})_{n},k(y))=0$$
pour tout $i>0$. Ceci implique donc que $(\tau_{\geq n}F'_{*})_{n}$
est localement libre dans un voisinage de $y$. Quitte \`a restreindre $Y$, on voit donc que 
le complexe $\tau_{\geq n}F'_{*}$ est parfait,
et de plus que le morphisme naturel
$$F'_{*} \longrightarrow \tau_{\geq n}F'_{*}$$
induit un quasi-isomorphisme 
$$F'_{*}\otimes^{\mathbb{L}}k(y) \longrightarrow \tau_{\geq n}F'_{*}\otimes^{\mathbb{L}}k(y).$$
Ceci implique alors que le c\^one $C$ du morphisme $F'_{*} \longrightarrow \tau_{\geq n}F'_{*}$
est tel que $C\otimes^{\mathbb{L}}k(y)$ soit contractile. Comme
$C$ est un complexe \`a cohomologie coh\'erente et born\'ee cela
implique par Nakayama que $C$ est quasi-isomorphe \`a $0$ dans un voisinage de $Y$. Ainsi, 
$F'_{*} \longrightarrow \tau_{\geq n}F'_{*}$ est un quasi-isomorphisme
dans un voisinage de $y$, ce qui implique que $F$ est parfait dans un 
voisinage de $y$. 
\hfill $\Box$ \\

Le lemme \ref{l7} appliqu\'e au complexe $\mathbb{R}f_{*}(E)$ 
termine alors la preuve de la proposition. \hfill $\Box$ \\

On d\'eduit de la proposition et de sa preuve la formule du changement
de bases suivante.

\begin{cor}\label{c1}
Soit 
$$\xymatrix{
X' \ar[d]_-{q} \ar[r]^-{g} & X \ar[d]^-{p} \\
Y' \ar[r]_-{f} & Y}$$
un diagramme cart\'esien d'espaces analytiques 
avec $p$ un morphisme propre et plat.
Alors, pour tout complexe parfait $E$ sur $X$, le morphisme naturel
$$\mathbb{L}f^{*}\mathbb{R}p_{*}(E) \longrightarrow \mathbb{R}q_{*}\mathbb{L}p^{*}(E)$$
est un isomorphisme dans $D(Y')$.
\end{cor}

\textit{Preuve:} Tout d'abord, le lemme \ref{l6} appliqu\'e \`a $F=k(y)$, le faisceau gratte-ciel
en un point $y\in Y$, implique que le corollaire \ref{c1} est vrai 
lorsque $Y'$ est un point. Les objets 
$\mathbb{L}f^{*}\mathbb{R}p_{*}(E)$ er $\mathbb{R}q_{*}\mathbb{L}p^{*}(E)$
\'etant parfaits sur $Y'$ (d'apr\`es la proposition \ref{p5}), le morphisme
$$\mathbb{L}f^{*}\mathbb{R}p_{*}(E) \longrightarrow \mathbb{R}q_{*}\mathbb{L}p^{*}(E)$$
est un isomorphisme dans $D(Y')$ si et seulement si pout tout $y'\in Y'$ d'inclusion $i_{y'} : \{y'\} \hookrightarrow Y'$, 
le morphisme induit
$$\mathbb{L}i_{y'}^{*}\mathbb{L}f^{*}\mathbb{R}p_{*}(E) \longrightarrow 
\mathbb{L}i_{y'}^{*}\mathbb{R}q_{*}\mathbb{L}p^{*}(E)$$
est un isomorphisme dans $D(\{y'\})$. Ainsi, les formules de changement de bases
pour $\{y'\} \longrightarrow Y'$ et $\{y'\} \longrightarrow Y'$ impliquent
le corollaire. \hfill $\Box$ \\

Nous terminons cette section par la d\'efinition de la version
dg des cat\'egories d\'eriv\'ees. Pour cela nous rappelons l'existence
de la cat\'egorie de mod\`eles $C(\mathcal{O}_{X})$  
des complexes de $\mathcal{O}_{X}$-modules sur un espace analytique $X$ (voir
l'appendice C). Cette cat\'egorie de mod\`eles est une
$C(\mathbb{C})$-cat\'egorie de mod\`eles au sens de \cite[4.2.18]{ho}, o\`u
$C(\mathbb{C})$ est la cat\'egorie de mod\`eles monoidale des
compelxes de $\mathbb{C}$-espaces vectoriels (pour la quelle les
\'equivalences sont les quasi-isomorphismes et les fibrations sont 
les \'epimorphismes). Nous consid\'erons alors
la dg-cat\'egorie $Int(C(\mathcal{O}_{X})$, des objets fibrants et cofibrants
dans $C(\mathcal{O}_{X})$ (voir \cite{to} pour plus de d\'etails sur la construction Int).

\begin{df}\label{ddg}
La \emph{dg-cat\'egorie d\'eriv\'ee des $\mathcal{O}_{X}$-modules} est 
$$L(X):=Int(C(\mathcal{O}_{X})).$$
La \emph{dg-cat\'egorie d\'eriv\'ee des $\mathcal{O}_{X}$-modules parfaits}
est la sous-dg-cat\'egorie pleine de $L(X)$ form\'ee des
complexes parfaits.
\end{df}

\section{Analytification des espaces de modules d'objets dans des dg-cat\'egories}

Dans cette section nous rappelons l'existence d'un espace
alg\'ebrique $M(B)$, classifiant les $B$-dg-modules \emph{simples} sur
une dg-alg\`ebre $B$ propre et lisse. Nous expliciterons aussi 
le foncteur des points de l'espace analytique $M(B)^{an}$. Nous supposerons
le lecteur famili\'e avec les notions de dg-alg\'ebres, dg-modules et 
de lissit\'e et propret\'e, tels que pr\'esent\'ees par exemple dans
\cite{tova}.

\subsection{Rappels sur les espaces de modules de dg-modules simples}

Consid\'erons une dg-alg\`ebre $B$ sur $\mathbb{C}$. On supposera que
$B$ est propre est lisse, c'est \`a dire qu'elle v\'erifie les deux conditions
suivantes (voir e.g. \cite{koso,tova}).

\begin{enumerate}
\item Le complexe sous-jacent \`a $B$ est parfait (i.e. $B$ est \`a cohomologie
born\'ee et de dimension finie).
\item $B$ est un objet compact dans $D(B\otimes B^{op})$, la cat\'egorie
d\'eriv\'ee des $B\otimes B^{op}$-dg-modules.
\end{enumerate}

Pour toute $\mathbb{C}$-alg\`ebre commutative $A$, non n\'ecessairement de type fini, 
on dispose d'une nouvelle $\mathbb{C}$-dg-alg\`ebre $B\otimes A$ et de sa cat\'egorie d\'eriv\'ee
$D(B\otimes A)$. Si $A\longrightarrow A'$ est un morphisme de $\mathbb{C}$-alg\`ebres
commutatives, le mophisme $B\otimes A \longrightarrow B\otimes A'$ induit un foncteur de changement de bases
$$-\otimes^{\mathbb{L}}_{A}A'=-\otimes^{\mathbb{L}}_{B\otimes A}B\otimes A' : D(B\otimes A) \longrightarrow
D(B\otimes A').$$
La cat\'egorie $B\otimes A-Mod$, des $B\otimes A$-dg-modules, poss\`ede un enrichissement naturel
dans la cat\'egorie mono\"\i dale $C(A)$ des complexes de $A$-modules. En effet, pour 
$E\in B\otimes A-Mod$ et $L\in C(A)$, on forme $L\otimes_{A}E$ qui est encore un 
$B\otimes A$-dg-module en utilisant la structure de $B$-dg-module sur $E$ et 
celle de $A$-dg-module sur $L$. En clair, l'action de $B\otimes A$ sur
$L\otimes_{A}E$ est donn\'ee par le morphisme
$$\xymatrix{
(B\otimes A)\otimes (L\otimes_{A}E) \simeq
(A\otimes L)\otimes_{A}(B\otimes E) \ar[r]^-{\mu_{L}\otimes \mu_{E}} &  L\otimes_{A}E,}$$
o\`u $\mu_{L}$ et $\mu_{E}$ sont les actions de $A$ et $B$ sur $L$ et $E$. Cet enrichissement 
de $B\otimes A-Mod$ dans $C(A)$ fait de $B\otimes A$ une
$C(A)$-cat\'egorie de mod\`eles au sens de \cite[4.2.18]{ho} (on peut aussi voir cela en 
remarquant que $B\otimes A-Mod$ est aussi la cat\'egorie
des $B\otimes A$-modules dans $C(A)$, lorsque $B\otimes A$ est vu 
comme un mono\"\i de dans $C(A)$ (i.e. comme une $A$-dg-alg\`ebre)). De ceci nous d\'eduisons
l'existence d'un enrichissement de $D(B\otimes A)$ dans $D(A)$. En particulier, on dispose de
$Hom$ \`a valeurs dans $D(A)$ que nous noterons
$$\mathbb{R}\underline{Hom}_{B}(-,-) : D(B\otimes A)^{op}\times D(B\otimes A) \longrightarrow D(A).$$

\begin{df}\label{d6}
\begin{enumerate}
\item Un objet $E\in D(B\otimes A)$ est \emph{parfait} si 
le complexe de $A$-modules sous-jacent \`a $E$ est un complexe parfait. 
\item un objet $E\in D(B\otimes A)$ est \emph{rigide}
si pour tout morphisme de $\mathbb{C}$-alg\`ebres commutatives $A\longrightarrow A'$,
et pour tout $i<0$ on a 
$$Ext^{i}(E\otimes^{\mathbb{L}}_{A}A',E\otimes^{\mathbb{L}}_{A}A')=0,$$
o\`u les $Ext$ sont calcul\'es dans $D(B\otimes A')$. 
\item Un objet $E\in D(B\otimes A)$ est \emph{simple} s'il est rigide et si de plus
pour tout morphisme de $\mathbb{C}$-alg\`ebres commutatives $A\longrightarrow A'$
le morphisme naturel
$$A' \longrightarrow Ext^{0}(E\otimes^{\mathbb{L}}_{A}A',E\otimes^{\mathbb{L}}_{A}A')$$
est un isomorphisme.
\end{enumerate}
\end{df}

\begin{rmk}\label{r5}
\emph{La notion d'objet parfait ci-dessus correspond en r\'ealit\'e \`a la notion
d'objet} pseudo-parfait \emph{de \cite{tova}. Cependant, comme $B$ est propre est lisse 
on sait que les objets pseudo-parfaits sont exactement les objets parfaits (voir \cite[Lem. 2.8 (3)]{tova}
pour une preuve de ce fait).}
\end{rmk}

Nous aurons besoin des quelques r\'esultats suivants.

\begin{prop}\label{p6}
Soient $A$ une $\mathbb{C}$-alg\`ebre commutative et \'ecrivons $A=Colim_{i}A_{i}$, o\`u la colimite est prise
sur l'ensemble filtrant des sous-$\mathbb{C}$-alg\`ebres de type fini $A_{i}\subset A$. 
\begin{enumerate}
\item Le foncteur naturel 
$$Colim (-\otimes^{\mathbb{L}}_{A}A_{i}) : Colim D(B\otimes A_{i}) \longrightarrow D(B\otimes A)$$
induit une \'equivalence sur les cat\'egories d'objets parfaits
$$Colim D_{parf}(B\otimes A_{i}) \simeq D_{parf}(B\otimes A).$$
\item Soit $E_{i}\in D(B\otimes A_{i})$ un objet parfait. Alors $E_{i}\otimes_{A_{i}}^{\mathbb{L}}A \in D(B\otimes A)$
est simple si et seulement s'il existe $A_{i} \subset A_{j} \subset A$, avec $A_{j}$ de type fini et 
telle que $E_{i}\otimes_{A_{i}}^{\mathbb{L}}A_{j} \in D(B\otimes A_{j})$ soit simple.
\end{enumerate}
\end{prop}

\textit{Preuve:} Le point $(1)$ est une cons\'equence de \cite[Lem. 2.10]{tova}. 
La suffisance de la condition dans $(2)$ 
se d\'eduit des d\'efinitions, car \^etre simple est clairement stable par changement de bases. Il nous
reste \`a montrer la n\'ecessit\'e. Pour cela, supposons que 
$E:=E_{i}\otimes_{A_{i}}^{\mathbb{L}}A$ soit simple dans $D(B\otimes A)$. On pose
$E_{j}:=E_{i}\otimes_{A_{i}}^{\mathbb{L}}A_{j}$ pour tout $j$ tel que $A_{i}\subset A_{j}$
(nous dirons alors que $j\geq i$). 

Pour tout $j\geq i$, on dispose d'un morphisme de complexes de $A_{j}$-modules
$$A_{j} \longrightarrow \mathbb{R}\underline{Hom}_{B}(E_{j},E_{j}),$$
o\`u le membre de droite fait r\'ef\'erence \`a l'enrichissement naturel 
de $D(B\otimes A_{j})$ dans $D(A_{j})$.  Comme le complexe de $A$-modules sous-jacent \`a $B\otimes A$ est parfait
et que $E_{j}$ est un $B\otimes A$-dg-module parfait on voit que le complexe de $A$-modules
$\mathbb{R}\underline{Hom}_{B}(E_{j},E_{j})$ est parfait. De plus, le fait que
$E_{j}$ soit parfait implique aussi que le morphisme naturel
$$\mathbb{R}\underline{Hom}_{B}(E_{i},E_{i}) \otimes_{A_{i}}^{\mathbb{L}}A_{j} \longrightarrow
\mathbb{R}\underline{Hom}_{B}(E_{j},E_{j})$$
est un quasi-isomorphisme. On d\'eduit de cela et du point $(1)$ que le morphisme naturel
$$\mathbb{R}\underline{Hom}_{B}(E_{i},E_{i}) \otimes_{A_{i}}^{\mathbb{L}}A \simeq
colim_{j\geq i} \mathbb{R}\underline{Hom}_{B}(E_{i},E_{i}) \otimes_{A_{i}}^{\mathbb{L}}A_{j} \longrightarrow
\mathbb{R}\underline{Hom}_{B}(E,E)$$
est un quasi-isomorphisme. 

On consid\`ere maintenant le morphisme naturel
$$A \longrightarrow \mathbb{R}\underline{Hom}_{B}(E,E),$$
qui est un morphisme de complexes parfaits de $A$-modules. On note
$K$ son c\^one. D'apr\`es ce que l'on a vu, on a
$K\simeq K_{i}\otimes^{\mathbb{L}}_{A_{i}}A$, o\`u $K_{i}$ est le c\^one du morphisme
$$A_{i} \longrightarrow \mathbb{R}\underline{Hom}_{B}(E_{i},E_{i}).$$
Par hypoth\`ese $E$ est simple, ce qui est \'equivalent au fait que
$K$ soit de Tor amplitude contenue dans $]0,\infty[$. 

\begin{lem}\label{l8}
Soit $A=colim_{i} A_{i}$ une colimite filtrante d'anneaux commutatifs et $K_{i}$ un complexe
parfait de $A_{i}$-modules pour un indice $i$ donn\'e. Si le complexe de $A$-modules
$K_{i}\otimes_{A_{i}}^{\mathbb{L}}A$ est de Tor amplitude contenue dans $[a,b]$
alors il existe un indice $j\geq i$ tel que $K_{i}\otimes_{A_{i}}^{\mathbb{L}}A_{j}$ soit de
Tor amplitude contenue dans $[a,b]$ en tant que complexe de $A_{j}$-modules
\end{lem}

\textit{Preuve du lemme \ref{l8}:} On supposera que $K_{i}$ est un complexe born\'e
de $A_{i}$-modules projectifs et de rangs finis. Par hypoth\`ese, il existe
un complexe de $A$-modules projectifs de type fini $L$, avec 
$L^{n}=0$ pour tout $n\notin [a,b]$, et une \'equivalence d'homotopie
$u : L \longrightarrow K$. Le complexe $L$ et
le morphisme $u$ sont alors tous deux d\'efinis sur
un $A_{j}$ avec $j\geq i$, et il existe donc un morphisme de complexes
$$u_{j} : L_{j} \longrightarrow K_{i}\otimes_{A_{i}}A_{j}$$
avec $u_{j}\otimes_{A_{j}}A=u$ (dans cette notation $L_{j}$ est un complexe de $A_{j}$-modules projectifs de type fini  avec 
$L_{j}^{n}=0$ pour tout $n\notin [a,b]$, et tel que $L_{j}\otimes_{A_{j}}A\simeq L$).
On peut faire de m\^eme avec un inverse \`a homotopie pr\`es $v : K \longrightarrow L$ de $u$. 
Quitte \`a prendre un $j$ plus grand on peut aussi trouver
$v_{j} : K_{i}\otimes_{A_{i}}A_{j} \longrightarrow L_{j}$ tel que $v_{j}\otimes_{A_{j}}A=v$. 
Enfin, les homotopies $h$ et $k$ reliant $vu$ et $uv$ avec les identit\'es sont elles aussi
d\'efinissables sur un $A_{j}$ pour $j$ assez grand. Ainsi, quitte \`a prendre $j$ assez grand
on voit que le morphisme 
$$u_{j} : L_{j} \longrightarrow K_{i}\otimes_{A_{i}}A_{j}$$
est une \'equivalence d'homotopie. Cela implique que
$K_{i}\otimes_{A_{i}}^{\mathbb{L}}A_{j}$ est de
Tor amplitude contenue dans $[a,b]$.
\hfill $\Box$ \\

Le lemme pr\'ec\'edent implique qu'il existe $A_{i} \subset A_{j} \subset A$ tel que
$K_{i}\otimes_{A_{i}}^{\mathbb{L}}A_{j}$ doit de Tor amplitude contenue dans
$]0,\infty[$. Comme nous avons vu que $K_{j}$ est le c\^one du morphisme
$$A_{j} \longrightarrow \mathbb{R}\underline{Hom}_{B}(E_{j},E_{j}),$$
nous en d\'eduisons que $E_{j}$ est simple. 
\hfill $\Box$ \\

\begin{cor}\label{c2}
Soient $A$ une $\mathbb{C}$-alg\`ebre commutative et \'ecrivons $A=Colim_{i}A_{i}$, o\`u la colimite est prise
sur l'ensemble filtrant des sous-$\mathbb{C}$-alg\`ebres de type fini $A_{i}\subset A$. Notons
$D_{pasi}(B\otimes A)$ et $D_{pasi}(B\otimes A_{i})$ les cat\'egories d'objets parfaits et simples. 
Alors le foncteur naturel
$$Colim D_{pasi}(B\otimes A_{i}) \longrightarrow D_{pasi}(B\otimes A)$$
est une \'equivalence.
\end{cor}

\textit{Preuve:} D\'ecoule de la proposition \ref{p6}. \hfill $\Box$ \\

\begin{prop}\label{p7}
Soit $(A,m)$ une $\mathbb{C}$-alg\`ebre locale et noeth\'erienne. Alors un objet
parfait $E\in D(B\otimes A)$ est simple si et seulement s'il satisfait aux deux conditions
suivantes.
\begin{enumerate}
\item Pour tout $i<0$ on a 
$$Ext^{i}(E\otimes^{\mathbb{L}}_{A}A/m,E\otimes^{\mathbb{L}}_{A}A/m)=0,$$
o\`u les $Ext$ sont calcul\'es dans $D(B\otimes A/m)$. 
\item Le morphisme naturel
$$A/m \longrightarrow Ext^{0}(E\otimes^{\mathbb{L}}_{A}A/m,E\otimes^{\mathbb{L}}_{A}A/m)$$
est un isomorphisme.
\end{enumerate}
\end{prop}

\textit{Preuve:} Comme nous l'avons utilis\'e lors de la preuve de la proposition \ref{p6}, 
$E$ est simple si et seulement si le c\^one du morphisme
$$A \longrightarrow \mathbb{R}\underline{Hom}_{B}(E,E)$$
est de Tor amplitude contenue dans $]0,\infty[$. Or, comme $A$ est noeth\'erien, on sait
qu'un $A$-module de type fini $M$ est plat si et seulement si 
$Tor_{n}^{A}(M,A/m)=0$ pour $n>0$. On d\'eduit facilement de cela qu'un complexe
parfait $K$ sur $A$ est de Tor amplitude contenue dans un 
intervalle $[a,b]$ si et seulement si $K\otimes_{A}^{\mathbb{L}}A/m$ est 
de Tor amplitude contenue dans $[a,b]$ (i.e. si et seulement si 
$Tor_{n}^{A}(K\otimes_{A}^{\mathbb{L}}A/m)=0$ pour $n\notin[-b,-a]$).
Ainsi, $E$ est simple si et seulement si $K$ est de Tor amplitude
contenue dans $]0,\infty[$, et donc si et seulement si 
$K\otimes_{A}^{\mathbb{L}}A/m$ est de Tor amplitude
contenue dans $]0,\infty[$. Comme $E$ est un $B\otimes A$-dg-module parfait on sait que
$K\otimes_{A}^{\mathbb{L}}A/m$  est le c\^one du morphisme induit
$$A/m \longrightarrow \mathbb{R}\underline{Hom}_{B}(E\otimes^{\mathbb{L}}_{A}A/m,E\otimes^{\mathbb{L}}_{A}A/m).$$
Ce qui termine la preuve de la proposition.
\hfill $\Box$ \\

Nous d\'efinissons maintenant un pr\'efaisceau $M'(B)$ sur $Aff$ de la fa\c{c}on suivante. 
Pour $A$ une $\mathbb{C}$-alg\`ebre de type finie, nous noterons
$M'(B)(A)$ l'ensemble des classes d'isomorphismes d'objets parfaits et simples
dans $D(B\otimes A)$. 
Les changements de bases pr\'eservent les objets parfaits et simples. Ainsi, 
un morphisme $A\longrightarrow A'$ entre $\mathbb{C}$-alg\`ebres de type fini induit une application
$$-\otimes^{\mathbb{L}}_{A}A' : M(B)(A) \longrightarrow M(B)(A').$$
Ceci d\'efinit un pr\'efaisceau $M'(B)\in Pr(Aff)$. 

\begin{df}\label{d7}
Le \emph{faisceau des $B$-dg-modules parfaits et simples} est le faisceau associ\'e
au pr\'efaisceau $M'(B)$. Il sera not\'e
$$M(B):=a(M'(B)).$$
\end{df}

Un important corollaire du th\'eor\`eme principal de \cite{tova} est la
repr\'esentabilit\'e de $M(B)$ par un espace alg\'ebrique localement de type fini.

\begin{thm}\label{t1}{(Voir \cite[Cor. 3.29, Rem. 3.30]{tova})}
Le faisceau $M(B)$ est un espace alg\'ebrique au sens de la d\'efinition \ref{d3}. 
Il est de plus quasi-s\'epar\'e (c'est \`a dire que le morphisme
diagonal $M(B) \longrightarrow M(B)\times M(B)$
est quasi-compact). 
\end{thm}

\subsection{Description de l'espace $M(B)^{an}$}

Comme pr\'ec\'edemment, soit $B$ une $\mathbb{C}$-dg-alg\`ebre
propre et lisse. 
Dans cette section nous nous proposons de donner une
description de l'espace analytique $M(B)^{an}$.

Rappelons que pour tout espace de Stein $S$, on dispose
d'une cat\'egorie de mod\`eles $C(\mathcal{O}_{S})$
des complexes de $\mathcal{O}_{S}$-modules (sur le petit
site des ouverts de Stein de $S$, voir l'appendice C). Cette cat\'egorie de mod\`eles
est de plus une $C(\mathbb{C})$-cat\'egorie de mod\`eles. On peut donc
parler de $B$-dg-modules dans $C(\mathcal{O}_{S})$: il s'agit d'objets
$E\in C(\mathcal{O}_{S})$ munis de morphismes
$$m : B\otimes E \longrightarrow E$$
qui sont unitaires et associatifs en un sens \'evident
(le produit tensoriel $B\otimes E$ est ici pris au-dessus de $\mathbb{C}$). Les morphismes
entre $B$-dg-modules dans $C(\mathcal{O}_{S})$ sont 
les morphismes dans $C(\mathcal{O}_{S})$ qui commutent avec les
morphismes structuraux $m$. Nous noterons la cat\'egorie des
$B$-dg-modules par $B-Mod(S)$. 

On peut aussi voir un $B$-dg-modules dans
$C(\mathcal{O}_{S})$ comme un module sur $B\otimes \mathcal{O}_{S}$ consid\'er\'e comme
mono\"\i de associatif et unitaire dans la cat\'egorie mono\"\i dale 
$(C(\mathcal{O}_{S}),\otimes_{\mathcal{O}_{S}})$. Ainsi, 
\cite{ss} implique l'existence d'une structure de cat\'egorie de mod\`eles
sur $B-Mod(S)$, dont les fibrations et les \'equivalences sont 
d\'efinies dans $C(\mathcal{O}_{S})$ par oubli de l'action de $B$. 
Nous noterons 
$$D(S,B):=Ho(B-Mod(S))$$
la cat\'egorie homotopique des $B$-dg-modules sur $S$. De plus, cette cat\'egorie
de mod\`eles est naturellement enrichie dans la cat\'egorie de mod\`eles
mono\"\i dale $C(\mathcal{O}_{S})$. On dispose ainsi de $Hom$ dans $D(S,B)$ \`a valeurs
dans $D(S) = Ho(C(\mathcal{O}_{S}))$ not\'es
$$\mathbb{R}\underline{Hom}_{B}(-,-) : D(S,B)^{op} \times D(S,B) \longrightarrow D(S).$$

Soit maintenant $f : S \longrightarrow S'$ un morphisme entre espaces 
de Stein. On dispose d'une adjonction de Quillen
$$f^{*} : C(\mathcal{O}_{S'}) \rightleftarrows C(\mathcal{O}_{S}) : f_{*}.$$
 Comme le foncteur $f^{*}$ est de plus muni d'une structure naturelle de 
foncteur mono\"\i dal sym\'etrique, cette adjonction induit une nouvelle
adjonction de Quillen
$$f^{*} : B-Mod(S') \rightleftarrows B-Mod(S) : f_{*}$$
sur les cat\'egories de mod\`eles de $B$-dg-modules.
On obtient donc une adjonction d\'eriv\'ee
$$\mathbb{L}f^{*} : D(S',B) \rightleftarrows D(S,B) : \mathbb{R}f_{*}.$$

\begin{df}\label{d8}
Soit $S$ un espace de Stein.
\begin{enumerate}
\item Un objet $E\in D(S,B)$ est \emph{parfait} si 
le complexe de $\mathcal{O}_{S}$-modules sous-jacent \`a $E$ est un complexe parfait. 
\item un objet $E\in D(S,B)$ est \emph{rigide}
si pour tout point $i_{s} : \{s\} \hookrightarrow S$, et pour tout $i<0$ on a 
$$Ext^{i}(\mathbb{L}i_{s}^{*}(E),\mathbb{L}i_{s}^{*}(E)):=0,$$
o\`u les groupes $Ext$ sont calcul\'es dans $D(\{s\},B)\simeq D(B)$. 
\item Un objet $E\in D(S,B)$ est \emph{simple} s'il est rigide et si de plus
pour tout point $i_{s} : \{s\} \hookrightarrow S$, le morphisme naturel
$$\mathbb{C} \longrightarrow Ext^{0}(\mathbb{L}i_{s}^{*}(E),\mathbb{L}i_{s}^{*}(E))$$
est un isomorphisme. 
\end{enumerate}
\end{df}

Nous d\'efinissons maintenant un foncteur
$$F_{B} : Ste^{op} \longrightarrow Ens$$
de la fa\c{c}on suivante. Pour $S\in Ste$, l'ensemble
$F_{B}(S)$ est le sous-ensemble des classes d'isomorphismes d'objets parfaits et simples
de $D(S,B)$.
Pour $f : S \longrightarrow S'$ un morphisme entre espaces de Stein, on 
dispose du changement de bases
$$\mathbb{L}f^{*} : D(S',B) \longrightarrow D(S,B).$$
Par d\'efinition ce foncteur pr\'eserve les conditions objets parfaits et simples et
donc induit
une application
$$\mathbb{L}f^{*} : F_{B}(S') \longrightarrow F_{B}(S).$$
Ceci d\'efinit le foncteur $F_{B}$. \\

\begin{prop}\label{p8}
Soit $S$ un espace de Stein et $E\in D(S,B)$ un objet parfait. 
Supposons qu'il existe un point $s\in S$ tel que 
$\mathbb{L}i_{s}^{*}(E)$ soit un objet simple dans $D(\{s\},B)\simeq D(B)$.
Alors il existe un voisinage ouvert de Stein $s\in U\subset S$ tel que 
la restriction de $E$ \`a $U$ soit un objet simple dans $D(U,B)$. 
\end{prop}

\textit{Preuve:} C'est le m\^eme principe que pour la proposition \ref{p7}, on montre que l'ensemble
des points $s$ de $S$ tels que $\mathbb{L}i_{s}^{*}(E)$ soit un objet simple est 
le lieu o\`u un certain complexe parfait est de Tor amplitude strictement positive. 
Pour cela, on consid\`ere la cat\'egorie de mod\`eles
$B-Mod(S)$. On rappelle que cette cat\'egorie est naturellement enrichie sur la cat\'egorie
mono\"\i dale $C(\mathcal{O}_{S})$ des complexes de $\mathcal{O}_{S}$-modules. 
Cet enrichissement fait de $B-Mod(S)$ une $C(\mathcal{O}_{S})$-cat\'egorie de mod\`eles. 
Les Hom enrichis de $B-Mod(S)$ \`a valeurs dans $C(\mathcal{O}_{S})$
seront not\'es $\underline{Hom}_{B}$, et leur version d\'eriv\'ee sera not\'ee
$\mathbb{R}\underline{Hom}_{B}$. 

Avec ces notations, on dispose d'un morphisme de complexes de $\mathcal{O}_{S}$-modules
$$\mathcal{O}_{S} \longrightarrow \mathbb{R}\underline{Hom}_{B}(E,E),$$
bien d\'efini dans $D(S)$. Nous commen\c{c}ons par remarquer que
$\mathbb{R}\underline{Hom}_{B}(E,E)$ est un complexe parfait sur $S$. En effet, 
on dispose d'un foncteur bi-exact entre cat\'egories triangul\'ees
$$D(B\otimes B^{op}) \times D(S,B) \longrightarrow D(S,B),$$
qui \`a $P\in D(B\otimes B^{op})$ et $F\in D(S,B)$ associe
$P\otimes^{\mathbb{L}}_{B}F$. Ce foncteur envoie 
le couple $(B,E)$ sur $E$. Or, comme $B$ est lisse, on sait que
$B$ appartient \`a la sous-cat\'egorie triangul\'ee et \'epaisse
engendr\'ee par $B\otimes B^{op}$. Cela implique que 
$E\simeq B\otimes^{\mathbb{L}}_{B}E$
appartient \`a la sous-cat\'egeorie triangul\'ee et \'epaisse
engendr\'ee par $(B\otimes B^{op})\otimes^{\mathbb{L}}_{B}E\simeq B\otimes E$.
On en d\'eduit que $E$ appartient \`a la sous-cat\'egorie 
triangul\'ee et \'epaisse engendr\'ee par les objets de la forme
$B\otimes F$, o\`u $F$ est un complexe parfait sur $S$. Ainsi, le complexe
$\mathbb{R}\underline{Hom}_{B}(E,E)$ appartient \`a la sous-cat\'egorie
triangul\'ee et \'epaisse de $D(S)$ engendr\'ee par 
les objets de la forme
$\mathbb{R}\underline{Hom}_{B}(B\otimes F,E)$. Or, comme
$B\otimes F$ est un $B$-dg-module libre, le complexe
$\mathbb{R}\underline{Hom}_{B}(B\otimes F,E)$ s'identifie naturellement au
complexe $\mathbb{R}\underline{Hom}(F,E)$ des morphismes de $F$ vers $E$ en tant que 
complexes de $\mathcal{O}_{S}$-modules. Comme $E$ et $F$ sont
parfaits, on en d\'eduit que 
$\mathbb{R}\underline{Hom}_{B}(E,E)$ est un complexe parfait. De plus, 
le m\^eme argument montre que pour tout point $s\in S$, le morphisme naturel
$$\mathbb{L}i_{s}^{*}(\mathbb{R}\underline{Hom}_{B}(E,E)) \longrightarrow
\mathbb{R}\underline{Hom}_{B}(\mathbb{L}i_{s}^{*}(E),\mathbb{L}i_{s}^{*}(E))$$
est un quasi-isomorphisme. 

Revenons au morphisme
$$\mathcal{O}_{S} \longrightarrow \mathbb{R}\underline{Hom}_{B}(E,E),$$
et notons $K$ son c\^one. C'est un complexe parfait, et par d\'efinition 
$\mathbb{L}i_{s}^{*}(E)$ est un objet simple si et seulement si
$H^{i}(\mathbb{L}i_{s}^{*}(K))=0$ pour tout $i\leq 0$. Cela est \'equivalent au fait
qu'au voisinage du point $s$ le complexe $K$ est quasi-isomorphe \`a un complexe
born\'e de fibr\'es vectoriels concentr\'e en degr\'es strictement positifs. Comme cela
est une condition ouverte, la proposition en d\'ecoule.
\hfill $\Box$ \\

Nous sommes maintenant en mesure de d\'emontrer la proposition suivante, qui 
donne une description de l'analytifi\'e du faisceau $M(B)$ des $B$-dg-modules simples.

\begin{prop}\label{p9}
Il existe un isomorphisme naturel entre $M(B)^{an}$ est le
faisceau associ\'e \`a $F_{B}$. 
\end{prop}

\textit{Preuve:} Notons $f : Aff \longrightarrow Ste$ le foncteur d'analytification et 
$f_{!} : Pr(Aff) \longrightarrow Pr(Ste)$ le foncteur induit sur les pr\'efaisceaux.
Rappelons que $M(B)$ est le faisceau associ\'e \`a un pr\'efaisceau
$M'(B)$. Il nous suffit donc de construire un morphisme de pr\'efaisceaux
$$\phi : f_{!}(M'(B)) \longrightarrow F_{B}$$
qui induise un isomorphisme sur les faisceaux associ\'es.

Pour construire $\phi$ il nous suffit par adjonction de construire 
un morphisme 
$$M'(B) \longrightarrow f^{*}(F_{B}).$$
Soit $X=Spec\, A \in Aff$. Nous disposons d'un foncteur
d'analytification 
$$-\otimes_{A}\mathcal{O}_{X^{an}} : C(A) \longrightarrow C(\mathcal{O}_{X^{an}}).$$
Ce foncteur est un foncteur de Quillen compatible avec l'enrichissement 
dans $C(\mathbb{C})$. Il induit donc un foncteur de Quillen \`a gauche
sur les cat\'egories de $B$-dg-modules
$$-\otimes_{A}\mathcal{O}_{X^{an}} : (A\otimes B)-Mod \longrightarrow B-Mod(X^{an}),$$
et donc un foncteur d\'eriv\'e
$$-\otimes_{A}^{\mathbb{L}}\mathcal{O}_{X^{an}} : 
D(A\otimes B) \longrightarrow D(X^{an},B).$$
Pour un point $x\in X(\mathbb{C})$, correspondant \`a un morphisme d'alg\`ebres
$A \longrightarrow \mathbb{C}$, le diagramme suivant
$$\xymatrix{
D(A\otimes B) \ar[rd]_-{-\otimes_{A}^{\mathbb{L}}\mathbb{C}} 
\ar[r]^-{-\otimes_{A}^{\mathbb{L}}\mathcal{O}_{X^{an}}} & D(X^{an},B) \ar[d]^-{\mathbb{L}i_{x}^{*}} \\
 & D(B)}$$
commute \`a isomorphisme pr\`es. Ceci implique que le foncteur
$-\otimes_{A}^{\mathbb{L}}\mathcal{O}_{X^{an}}$ envoie $M'(B)(A)$ dans 
$F_{B}(X^{an})$, et donc induit une application
$$\phi : M'(B)(A) \longrightarrow f^{*}(F_{B})(A).$$
Cette application est fonctorielle en $A$ et donc fournit un morphisme
de pr\'efaisceaux
$$\phi : M'(B) \longrightarrow f^{*}(F_{B})$$
et par adjonction un morphisme de pr\'efaisceaux sur $Ste$
$$\phi : f_{!}(M'(B)) \longrightarrow F_{B}.$$
Il nous reste \`a voir que ce morphisme induit un isomorphisme fibres \`a fibres. 
Pour cela, soit $X\in Ste$, $x\in X$, et consid\'erons le morphisme induit
$$f_{!}(M'(B))(X_{x}) \longrightarrow F_{B}(X_{x}).$$
D'apr\`es la proposition \ref{p4} ce morphisme est isomorphe \`a
$$Colim_{A\subset \mathcal{O}_{x}}M'(B)(A) \longrightarrow 
Colim_{x\in U\subset X}F_{B}(U).$$

Le corollaire \ref{c2} permet de calculer la colimite de gauche. En effet, elle
implique que l'on a
$$Colim_{A\subset \mathcal{O}_{x}}M'(B)(A)\simeq M'(B)(\mathcal{O}_{x}).$$
Il nous reste donc \`a montrer que le morphisme naturel
$$M'(B)(\mathcal{O}_{x}) \longrightarrow F_{B}(X_{x})$$
est une bijection. Nous allons commencer par d\'ecrire ce morphisme. 
Pour cela, soit $x\in U\subset X$ un voisinage ouvert de Stein de $x$ dans $X$. 
Soit $L(U):=Int(C(\mathcal{O}_{U}))$ la dg-cat\'egorie des objets fibrants et 
cofibrants dans $C(\mathcal{O}_{U})$ (voir \cite{to} pour plus de d\'etails
sur la construction $Int$), et $L_{parf}(U)$ la sous-dg-cat\'egorie
pleine de $L(U)$ form\'ee des complexes parfaits
de $\mathcal{O}_{U}$-modules. D'apr\`es \cite[Thm. 4.2]{to} nous avons une bijection naturelle
entre les classes d'isomorphismes d'objets dans $D(U,B)$ et l'ensemble
des morphismes $[B,L(U)]$ dans la cat\'egorie homotopique
des dg-cat\'egories. Ainsi, $F_{B}(U)$ s'identifie naturellement
au sous-ensemble de $[B,L_{parf}(U)]$ form\'e des morphismes
$B\longrightarrow L_{parf}(U)$ tels que pour tout point $s\in S$, 
le $B$-dg-modules correspondant au morphisme
$$\xymatrix{B \ar[r] & L_{parf}(U) \ar[r]^-{\mathbb{L}i_{s}^{*}} & L_{parf}(\{s\})}$$
soit un objet simple dans $D(B)$.
Cela implique donc que $F_{B}(X_{x})$ s'identifie naturellement \`a un sous-ensemble
de $Colim_{x\in U\subset X}[B,L_{parf}(U)]$. Or, $B$ \'etant 
propre et lisse, c'est une dg-cat\'egorie homotopiquement de pr\'esentation finie
(voir \cite[Cor. 2.13]{tova}), et donc nous avons une bijection
$$Colim_{x\in U\subset X}[B,L_{parf}(U)]\simeq [B,Colim_{x\in U\subset X}L_{parf}(U)].$$
Consid\'erons maintenant $U$ un voisinage ouvert de Stein de $x$ et le foncteur fibre en $x$
$$C(\mathcal{O}_{U}) \longrightarrow C(\mathcal{O}_{x}).$$
Ce foncteur est de Quillen \`a gauche, et comme tous les objets dans 
$C(\mathcal{O}_{x})$ sont fibrants il induit un morphisme de dg-cat\'egories
$$Int(C(\mathcal{O}_{U})) \longrightarrow Int(C(\mathcal{O}_{x})).$$
En se retreignant aux objets parfaits on trouve un morphisme
$$L_{parf}(U) \longrightarrow L_{parf}(\mathcal{O}_{x})$$
des complexes parfaits de $\mathcal{O}_{U}$-modules vers les complexes
parfaits sur l'anneau $\mathcal{O}_{x}$. Ce morphisme induit un morphisme de dg-cat\'egories
$$Colim_{x\in U\subset X}L_{parf}(U) \longrightarrow L_{parf}(\mathcal{O}_{x})$$
qui est une quasi-\'equivalence. Ainsi, d'apr\`es ce que l'on a vu, 
$F_{B}(X_{x})$ s'identifie donc \`a un sous-ensemble de
$$[B,Colim_{x\in U\subset X}L_{parf}(U)]\simeq [B,L_{parf}(\mathcal{O}_{x})].$$
D'apr\`es la proposition \ref{p8} il n'est pas difficile de voir que 
ce sous-ensemble est pr\'ecis\`ement le sous-ensemble
des classes d'isomorphismes de $D(B\otimes \mathcal{O}_{x})$ form\'e
des objets compacts $E$ qui sont tels que 
$E\otimes_{\mathcal{O}_{x}}^{\mathbb{L}}\mathcal{O}_{x}/m_x$
soit un objet simple dans $D(B)$. D'apr\`es la proposition \ref{p7} ce sous-ensemble correspond 
aussi au sous-ensemble $M'(B)(\mathcal{O}_{x})$ de $[B,L_{parf}(\mathcal{O}_{x})]$. 
Ceci termine la preuve que l'application
$$M'(B)(\mathcal{O}_{x}) \longrightarrow F_{B}(X_{x})$$
est bijective, et donc termine la preuve de la proposition.
\hfill $\Box$ \\

\section{Le th\'eor\`eme de caract\'erisation}

Nous arrivons maintenant \`a l'\'enonc\'e du th\'eor\`eme principal de ce travail. 
Pour cela, rappelons qu'une dg-cat\'egorie $T$ est \emph{satur\'ee} s'il existe une
dg-alg\`ebre propre et lisse $B$ et une quasi-\'equivalence $T\simeq L_{parf}(B)$
entre $T$ est la dg-cat\'egorie des $B$-dg-modules cofibrants et parfaits 
(voir \cite{tova} pour plus de d\'etails sur
cette notion, o\`u cette dg-cat\'egorie est not\'e $\widehat{B^{op}}_{pe}$). 

\begin{thm}\label{t2}
Soit $X$ un espace analytique connexe, compact et lisse. Alors $X$ est alg\'ebrisable si et seulement si
la dg-cat\'egorie $L_{parf}(X)$ est satur\'ee.
\end{thm}

La n\'ecessit\'e est bien connue. Lorsque $X$ est un sch\'ema propre et lisse cela est d\'emontr\'e
dans \cite[Lem. 3.27]{tova} (noter que l'on a $L_{parf}(X)\simeq L_{parf}(X^{an})$ d'ap\`es GAGA). 
Dans le cas g\'en\'eral o\`u $X$ est un espace alg\'ebrique propre et lisse on peut 
utiliser le lemme de Chow pour montrer que $L_{parf}(X)$ est satur\'ee (voir
l'appendice B pour les d\'etails). \\

Le reste de cette section est consacr\'e \`a la preuve de la suffisance de l'\'enonc\'e du th\'eor\`eme.
Pour cela on se fixe $X$ un espace analytique compact et lisse tel que $L_{parf}(X)$
soit satur\'ee. On se fixe aussi une quasi-\'equivalence $L_{parf}(B)\simeq L_{parf}(X)$ avec
$B$ une dg-alg\`ebre propre et lisse. Cela revient \`a se fixer un 
objet $E\in D_{parf}(X)$, dont l'enveloppe triangul\'ee et \'epaisse
est $D_{parf}(X)$, et tel que $\mathbb{R}\underline{End}(E)$ soit
propre et lisse. Nous nous fixerons plus pr\'ecis\'emment 
un objet fibrant et cofibrant $E\in C(\mathcal{O}_{X})$ dans la cat\'egorie
de mod\`eles des complexes de $\mathcal{O}_{X}$-modules qui soit un 
repr\'esentant pour $E$. Nous prendrons alors $B:=\underline{End}(E)$
la dg-alg\`ebre des endomorphismes de l'objet $E$. 

\subsection{Alg\'ebrisation de $M(X)$}

Nous allons d\'efinir un faisceau $M(X)$, associ\'e \`a un pr\'efaisceau
$M'(X)$. 

\begin{df}\label{d9}
Soit $S\in Ste$. Nous dirons que   
$E\in D_{parf}(X\times S)$ est \emph{simple} s'il v\'erifie les deux conditions suivantes:
\begin{enumerate}
\item pour tout $s\in S$, le morphisme naturel 
$$\mathbb{C} \longrightarrow Ext^{0}(\mathbb{L}j_{s}^{*}(E),\mathbb{L}j_{s}^{*}(E))$$
est un isomorphisme, 
\item pour tout $s\in S$ et tout $i<0$, on a 
$$Ext^{i}(\mathbb{L}j_{s}^{*}(E),\mathbb{L}j_{s}^{*}(E))=0,$$
o\`u $j_{s} : X\times \{s\} \hookrightarrow X\times S$ est l'inclusion
naturelle et les groupes $Ext^{i}$ sont calcul\'es dans $D(X)$.
\end{enumerate}
\end{df}

On d\'efinit alors $M'(X)(S)$ comme \'etant l'ensemble
des classes d'isomorphismes d'objets simples dans
$D_{parf}(X\times S)$. 
Pour un morphisme $f : S \longrightarrow S'$ le changement de bases
$$\mathbb{L}f^{*} : D_{parf}(X\times S') \longrightarrow D_{parf}(X\times S)$$
pr\'eserve les deux propri\'et\'es pr\'ec\'edentes et induit donc une application
$$\mathbb{L}f^{*} : M'(X)(S') \longrightarrow M'(X)(S).$$
Ceci d\'efinit un pr\'efaisceau $M'(X)$ sur $Ste$.

\begin{df}\label{d10}
Le \emph{faisceau des complexes parfaits simples sur $X$} est 
le faisceau associ\'e au pr\'efaisceau $M'(X)$. Il sera not\'e 
$M(X)$.
\end{df}

Le but de cette partie est de d\'emontrer la proposition suivante.

\begin{prop}\label{p10}
Le faisceau $M(X)$ est un espace analytique alg\'ebrisable. 
\end{prop}

\textit{Preuve:} Nous allons montrer plus pr\'ecisemment qu'il existe
un isomorphisme de faisceaux
$$M(X)\simeq M(B)^{an}.$$
Comme nous savons que $M(B)$ est un espace alg\'ebrique (voir \ref{t1}) cela impliquera la
proposition.

Pour cela nous \'etudierons la situation g\'en\'erale suivante. 
Soit $\mathcal{C}$ une cat\'egorie de mod\`eles sym\'etrique mono\"\i dale
qui est stable et qui v\'erifie les conditions de \cite{ss}. Soit 
$M$ une $\mathcal{C}$-cat\'egorie de mod\`eles, 
dont nous noterons $\underline{Hom}$ les morphismes enrichis
dans $\mathcal{C}$. Soit $E_{0}\in M$ un objet cofibrant et 
notons $B_{0}:=\underline{End}(E_{0})$ le mono\"\i de des endomorphismes de $E_{0}$ dans $\mathcal{C}$.
D'apr\`es les r\'esultats de \cite{ss}, il existe sur $B_{0}-Mod$, la cat\'egorie des $B_{0}$-modules
dans $\mathcal{C}$, une structure de cat\'egorie de mod\`eles pour
la quelle les fibrations et les \'equivalences sont d\'efinies dans $\mathcal{C}$. 
On d\'efinit alors un foncteur
$$\phi_{E_{0}} : M \longrightarrow B_{0}-Mod$$
qui associe \`a un objet $X\in M$ l'objet $\underline{Hom}(E_{0},X)$
muni de l'action naturelle \`a gauche du mono\"\i de
$\underline{End}(E_{0})$. Ce foncteur est un foncteur de Quillen \`a droite, dont 
l'adjoint \`a gauche sera not\'e $\psi_{E_{0}}$. 
Nous en d\'eduisons une adjonction 
$$\mathbb{L}\psi_{E_{0}} : Ho(B_{0}-Mod) \rightleftarrows Ho(M) : \mathbb{R}\phi_{E_{0}}.$$

Soit maintenant $S\in Ste$ un espace de Stein. Nous appliquons
les consid\'erations pr\'ec\'edentes \`a la $C(\mathcal{O}_{S})$-cat\'egorie de mod\`eles
$C(\mathcal{O}_{X\times S})$, et \`a $E_{0}:=E\boxtimes \mathcal{O}_{S}$
(rappelons que $E\in C(\mathcal{O}_{X})$
est un repr\'esentant fibrant et cofibrant d'un g\'en\'erateur
de $L_{parf}(X)$). Nous avons donc un foncteur 
$$\phi_{E_{0}} : D(X\times S) \longrightarrow D(B_{0}-Mod),$$
o\`u $B_{0}=\underline{End}(E_{0})$. On dispose de plus d'un 
morphisme naturel de mono\"\i des dans $C(\mathcal{O}_{S})$
$$B\otimes \mathcal{O}_{Z}=\underline{End}(E)\otimes \mathcal{O}_{Z} \longrightarrow \underline{End}(E_{0})$$
induit par le foncteur image inverse
$$C(\mathbb{C}) \longrightarrow C(\mathcal{O}_{S})$$
associ\'e au morphisme $S \longrightarrow *$. Ce morphisme induit un foncteur d'oubli sur
les cat\'egories de modules
$$D(B_{0}-Mod) \longrightarrow D(B-Mod)=D(S,B).$$
En composant avec $\phi_{E_{0}}$ on obtient un foncteur
$$\phi : D(X\times S) \longrightarrow D(S,B).$$
Ce foncteur est le compos\'e de deux foncteurs adjoints \`a droite, et donc poss\`ede un adjoint \`a gauche
$$\psi : D(S,B) \longrightarrow D(B_{0}-Mod).$$

Pour $F\in D(X\times S)$, l'objet de $D(S)$ sous-jacent \`a 
$\phi(F)$ est $\mathbb{R}p_{*}(E_{0}^{\vee}\otimes^{\mathbb{L}} F)$, o\`u 
$E_{0}^{\vee}$ est le dual du complexe parfait $E_{0}$ et 
$p : X\times S \longrightarrow S$ est la projection sur le second facteur.
Comme $E_{0}^{\vee}$ est parfait sur $X\times S$, la proposition \ref{p5} implique que le foncteur
$\phi$ envoie $D_{parf}(X\times S)$ dans la 
sous-cat\'egorie pleine de $D(S,B)$ form\'ee des objets qui sont 
parfaits comme complexes de $\mathcal{O}_{S}$-modules.

Nous avons donc un foncteur induit entre cat\'egories d'objets parfaits
$$\phi : D_{parf}(X\times S) \longrightarrow D_{parf}(S,B).$$
La propri\'et\'e de changement de bases (voir corollaire \ref{c1}) implique que 
$\phi$ est compatible, \`a isomorphisme pr\`es, aux changement de bases
induit par des morphismes $S' \longrightarrow S$. De plus, 
comme $E$ est g\'en\'erateur de $L_{parf}(X)$, le foncteur
$$\phi : D_{parf}(X\times S) \longrightarrow D_{parf}(S,B)$$
est une \'equivalence lorsque $S=*$ est r\'eduit \`a un point. Ces deux propri\'et\'es
impliquent que le foncteur
$$\phi : D_{parf}(X\times S) \longrightarrow D_{parf}(S,B)$$
pr\'eserve aussi les objets simples au sens des d\'efinitions
\ref{d8} et \ref{d9}. Il induit donc un 
morphisme de faisceaux sur $Ste$
$$\phi : M(X) \longrightarrow M(B)^{an}.$$

Il nous reste \`a montrer que $\phi$ est un isomorphisme. Pour cela, 
nous allons montrer que pour tout $S\in Ste$, le foncteur
$$\phi : D_{parf}(X\times S) \longrightarrow D_{parf}(S,B)$$
est une \'equivalence. On consid\`ere le foncteur
$$\psi : D(S,B) \longrightarrow D(X\times S),$$
adjoint \`a gauche du foncteur $\phi : D(X\times S) \longrightarrow D(S,B)$.

\begin{lem}\label{l9}
Le foncteur
$$\psi : D(S,B) \longrightarrow D(X\times S)$$
pr\'eserve les objets parfaits.
\end{lem}

\textit{Preuve du lemme:} Comme nous l'avons vu lors de la preuve de la proposition \ref{p8},  
$D_{parf}(S,B)$ est la plus petite sous-cat\'egorie triangul\'ee
\'epaisse de $D(S,B)$ contenant les objets de la forme
$B\otimes F$, o\`u $F\in D_{parf}(S)$ (avec la structure de
$B$-dg-module est donn\'ee par l'action de $B$ sur elle-m\^eme). 
Ainsi, pour montrer que $\psi$ pr\'eserve les
objets parfaits il suffit de montrer que 
$\psi(B\otimes F)$ est parfait
sur $X\times S$ lorsque $F$ est parfait sur $S$. Or, on voit par adjonction
que l'on a 
$$\psi(B\otimes F)\simeq E\boxtimes^{\mathbb{L}} F.$$
\hfill $\Box$ \\

D'apr\`es le lemme on dispose d'une adjonction
$$\psi : D_{parf}(S,B) \rightleftarrows D_{parf}(X\times S) : \phi.$$
De plus, cette adjonction est compatible aux changements de bases sur $S$. Ainsi, 
pour voir que les morphismes d'adjonction
$$Id \longrightarrow \phi \psi \qquad
\psi\phi \longrightarrow Id$$
sont des isomorphismes on peut supposer que $S$ est \'egal \`a un point
(cela est justifiable du fait que tous les objets consid\'er\'es sont des complexes parfaits
sur $S$ et sur $X\times S$, et donc qu'\^etre un quasi-isomorphisme
peut se tester points par points). Mais dans ce cas, le fait que 
cette adjonction soit une \'equivalence d\'ecoule  du fait que
$E$ soit un g\'en\'erateur de la dg-cat\'egorie $L_{parf}(X)$. 

Nous venons de voir que 
$$\phi : D_{parf}(X\times S) \longrightarrow D_{parf}(S,B)$$
\'etait une \'equivalence compatible aux changements de bases en $S$. 
Cela implique qu'en passant aux classes d'isomorphismes d'objets simples
ce foncteur induit un isomorphisme de faisceaux
$$\phi : M(X) \simeq M(B)^{an}.$$
\hfill $\Box$ \\

\subsection{Plongement de $X$ dans $M(X)$}

Le but de cette section est de d\'efinir un plongement ouvert 
de $X$ dans l'espace alg\'ebrisable $M(X)\simeq M(B)^{an}$. 
Pour cela on envoie un point $x$ de $X$ dans le 
faisceau gratte-ciel $k(x)\in D_{parf}(X)$, qui est un objet simple
et donc un point de $M(X)$. 

\begin{prop}\label{p11}
Il existe un morphisme injectif et \'etale
$$j : X \longrightarrow M(X).$$
\end{prop}

\textit{Preuve:} Soit $S\in Ste$ un espace de Stein et $u : S \longrightarrow X$ 
un morphisme, consid\'er\'e comme un \'el\'ement de $X(S)$. Notons
$\Gamma(u) \subset X\times S$ le  graphe de $u$, qui est un 
sous-ensemble analytique de $X\times S$ dont nous noterons
$I(u)\subset \mathcal{O}_{X\times S}$ l'id\'eal. Notons
$F(u)$ le faisceau coh\'erent $\mathcal{O}_{X\times S}/I(u)$, que l'on consid\`ere
comme un objet dans $D_{parf}(X\times S)$. L'objet 
$F(u)$ est simple au sens de la d\'efinition \ref{d9}, car pour
$s\in S$ on a
$$\mathbb{R}j_{s}^{*}(F(u))\simeq k(u(s)),$$
o\`u $k(u(s))$ est le faisceau gratte-ciel centr\'e au point $u(s)\in X$. 
La construction $u \mapsto F(u)$ induit ainsi une application
$$j : X(S) \longrightarrow M'(X)(S).$$
Lorsque $S$ varie dans $Ste$ cela d\'efinit un morphisme de pr\'efaisceaux
$$j : X \longrightarrow M'(X)$$
et donc de faisceaux
$$j : X\longrightarrow M(X).$$

La morphisme $j$ est injectif car si deux points
$x$ et $x'$ dans $X$ sont tels que $k(x)$ et $k(x')$ soient
isomorphes dans $D(X)$ alors $x=x'$ pour des raisons de supports.  Il nous reste \`a montrer que
$j$ est \'etale. 

Nous commencerons par montrer que $j$ est un monomorphisme. Pour cela, 
soit $u,v : S \longrightarrow X$ deux morphismes avec $S\in Ste$ tels que
$j(u)=j(v)$.
Dire que $j(u)=j(v)$ dans $M(X)(S)$ est \'equivalent \`a dire que 
les faisceaux structuraux des graphes $\Gamma(u)$ et $\Gamma(v)$ de $u$ et $v$ sont 
isomorphes comme faisceaux coh\'erents sur $X\times S$, au moins
localement sur $S$. Comme la propri\'et\'e $u=v$ est locale sur
$S$, on peut supposer que faisceaux structuraux des graphes $\Gamma(u)$ et $\Gamma(v)$ sont 
isomorphes. Cela implique facilement que les sous-espaces analytiques d\'efinis
par ces graphes sont \'egaux (en tant que sous-espaces analytiques 
de $X\times S$). Ceci implique bien entendu que les morphismes $u$ et $v$ sont
\'egaux. Ainsi, $j$ est non-seulement injectif mais aussi un monomorphisme.
Il nous reste donc \`a montrer que $j$ est aussi formellement lisse, ce qui 
impliquera qu'il est \'etale.

Soit donc $S\in Ste$ un espace de Stein et $i : S_{0}\subset S$ un sous-espace
analytique ferm\'e d\'efini par un id\'eal de carr\'e nul. Supposons que l'on ait un 
diagramme commutatif de faisceaux
$$\xymatrix{
X \ar[r]^-{j} & M(X) \\
S_{0} \ar[u]^-{u} \ar[r]_-{i} & S. \ar[u]_-{v}}$$
On cherche \`a montrer qu'il existe un morphisme
$w : S \longrightarrow X$ faisant commuter le diagramme
$$\xymatrix{
X \ar[r]^-{j} & M(X) \\
S_{0} \ar[u]^-{u} \ar[r]_-{i} & \ar[ul]^-{w}S \ar[u]_-{v}}$$
Comme nous savons d\'ej\`a que $j$ est un monomorphisme, et donc
non-ramifi\'e, l'existence de ce morphisme $w$ est une condition locale sur
$S$. On peut donc supposer que le morphisme
$v : S \longrightarrow M(X)$ se factorise par un morphisme
$v' : S \longrightarrow M'(X)$ de sorte \`a ce que le diagramme
$$\xymatrix{
X \ar[r]^-{j} & M'(X) \\
S_{0} \ar[u]^-{u} \ar[r]_-{i} & S \ar[u]_-{v'}}$$
soit un diagramme commutatif de pr\'efaisceaux. La donn\'ee de ce dernier diagramme
est \'equivalente \`a la donn\'ee d'un couple $(u,E)$, o\`u 
\begin{enumerate}
\item $u : S_{0} \longrightarrow X$ est un morphisme dont le faisceau
structural du graphe $\Gamma(u)\subset X\times S_{0}$ sera not\'e
$E_{0}(u)$.
\item $E \in D_{parf}(X\times S)$ est un objet simple tel que
$\mathbb{L}i^{*}(E)$ et $E_{0}(u)$ soit isomorphes dans $D_{parf}(X\times S_{0})$.
\end{enumerate}
Pour montrer que le rel\`evement $w$ comme ci-dessus existe il suffit 
de montrer qu'il existe un morphisme $w' : S \longrightarrow X$, dont le
faisceau structural du graphe sera not\'e $E(w')\in D_{parf}(X\times S)$, 
et qui est tel que
$$E(w')\simeq E \qquad w'\circ i=u.$$
Comme le morphisme $j$ est un monomorphisme il nous suffit en r\'ealit\'e de
montrer l'existence de $w$' tel que $E(w')\simeq E$, la seconde condition \'etant alors automatique.

Pour commencer, nous avons $E\in D_{parf}(X\times S)$ tel que
$\mathbb{L}i^{*}(E)$ soit un faisceau coh\'erent et plat sur
$S_{0}$. Comme dans \cite[Lem. 5.2]{anto} cela implique que $E$ est lui-m\^eme un faisceau
coh\'erent et plat sur $S$. On consid\`ere alors la projection
$p : X\times S \longrightarrow S$, et l'image directe
$\mathbb{R}p_{*}(E)$ de $E$ sur $S$. Par la propri\'et\'e de changement 
de bases du corollaire \ref{c1} on voit que pour tout $s\in S$, le complexe
$\mathbb{L}j_{s}^{*}\mathbb{R}p_{*}(E)$ est quasi-isomorphe \`a $ \mathbb{C}$ concentr\'e
en degr\'e z\'ero. Cela implique que 
$\mathbb{R}p_{*}(E)$ est un fibr\'e en droite $L$ sur $S$. Quitte \`a restreindre $S$ si n\'ecessaire
on peut aussi supposer que $L$ est trivial et choisir une trivialisation
$e : \mathcal{O}_{S} \simeq \mathbb{R}p_{*}(E)$. Le morphisme $e$ correspond par adjonction \`a un 
morphisme de faisceaux coh\'erents $f : \mathcal{O}_{X\times S} \longrightarrow E$.
En utilisant que $\mathbb{L}i^{*}(E)\simeq E_{0}(u)$ on voit, points par points sur $X\times S$, que
le morphisme $f$ est un \'epimorphisme. Ainsi, $E$ est-il isomorphe au faisceau 
structural d'un sous-espace analytique $Z\subset X\times S$, plat sur $S$. 
Enfin, le morphisme $Z \longrightarrow S$ est plat et tel que le morphisme induit
$$Z\times_{S}S_{0}\simeq \Gamma(u) \longrightarrow S_{0}$$
soit un isomorphisme, et donc est lui-m\^eme un isomorphisme. Cela implique l'existence
d'un morphisme $w' : S \longrightarrow X$ dont le graphe est \'egal au sous-espace
$Z$. On a donc bien $E(w')\simeq \mathcal{O}_{Z}\simeq E$. 
\hfill $\Box$ \\

\subsection{Alg\'ebrisation de $X$}

Pour achever la preuve du th\'eor\`eme il nous reste \`a d\'emontrer la proprosition suivante.

\begin{prop}\label{p12}
Soit $M$ un espace alg\'ebrique quasi-s\'epar\'e
(i.e. le morphisme diagonal $M\longrightarrow M\times M$ est 
quasi-compact). 
S'il existe un morphisme injectif et \'etale
$$j : X \longrightarrow M^{an}$$
alors $X$ est alg\'ebrisable (i.e. il existe un espace
alg\'ebrique $Y$ tel que $Y^{an}\simeq X$).
\end{prop}

\textit{Preuve:} Soit $\{U_{i} \longrightarrow M\}$ une famille couvrante
de morphismes \'etales avec $U_{i}$ des sch\'emas affines. Comme
$M$ est quasi-s\'epar\'e chacun des morphismes
$U_{i} \longrightarrow M$ est \'etale et de type fini. L'image de
$U_{i} \longrightarrow M$ est donc un sous-espace 
alg\'ebrique ouvert $M_{i}\subset M$ de type fini. 
Ainsi, l'espace alg\'ebrique $M$ est une r\'eunion de sous-espaces
ouverts de type fini. Par compacit\'e de $X$ on peut donc supposer que 
$j$ se factorise en $X \longrightarrow (M')^{an} \subset M^{an}$, avec $M'$ un sous-espace
alg\'ebrique ouvert de type fini. En d'autres termes on peut supposer que
$M$ est de type fini.

L'espace analytique $X$ est irr\'eductible (au sens analytique), et 
donc le morphisme $j : X \longrightarrow M^{an}$ se factorise
en $X \longrightarrow Z \subset M^{an}$, o\`u $Z$
est une composante irr\'eductible analytique de $M^{an}$. 
Comme le morphisme compos\'e 
$$X \longrightarrow Z \hookrightarrow M^{an}$$
est \'etale et que $Z \longrightarrow M^{an}$ est un monomorphisme, on voit 
que le morphisme $X \longrightarrow Z$ est encore injectif et \'etale. De plus, 
d'apr\`es \cite[Exp. XII, Prop. 2.4]{sga1}, il existe
une composante irr\'eductible $M_{0}$ de $M$ (au sens alg\'ebrique)
tel que $Z\simeq M_{0}^{an}$. Cela implique que l'on peut supposer
que $M$ est irr\'eductible. De la m\^eme fa\c{c}on, comme $X$ est r\'eduit on
pourra supposer que $M$ est r\'eduit.

L'immersion ouverte $j : X \hookrightarrow M^{an}$
induit un morphisme injectif sur les corps de fonctions m\'eromorphes
$$K(M^{an}) \hookrightarrow K(X).$$
De plus, comme $j$ est une immersion ouverte nous avons les \'egalit\'es 
$$dim(X) = dim(M^{an}) = dim(M)$$ 
(voir \cite[Exp. XII, Prop. 2.1]{sga1}). 
Comme $M$ est quasi-s\'epar\'e il existe un ouvert affine dense 
$U$ de $M$ (\cite[Prop. I.5.19]{Kn}) qui est de dimension $n=dim(X)$. Ainsi, on trouve un 
plongement
$$K(U)\simeq K(M) \subset K(M^{an}) \hookrightarrow K(X),$$
o\`u $K(U)$ et $K(M)$ d\'esigne les corps de fractions
rationnelles sur $U$ et $M$. 
Comme $K(U)$ est de degr\'e de transcendance \'egal \`a $n$, l'espace
analytique $X$ est un espace de Moishezon, et donc est alg\'ebrisable
d'apr\`es \cite{ar}.
\hfill $\Box$ \\

Pour finir la preuve du th\'eor\`eme \ref{t2} nous appliquons
la proposition \ref{p12} au morphisme injectif et \'etale
$$j : X \longrightarrow M(X)\simeq M(B)^{an}$$
de la proposition \ref{p11}. Ceci est possible car $M(B)$ est un 
espace alg\'ebrique  
quasi-s\'epar\'e d'apr\`es le th\'eor\`eme \ref{t1}.

\begin{appendix}

\section{Dg-cat\'egories satur\'ees et cat\'egories triangul\'ees satur\'ees}

Dans cette section nous comparons les notions de saturation dans les cadres
des dg-cat\'egories et des cat\'egories triangul\'ees, et nous comparons ainsi les
notions de \cite{bv} et \cite{tova}. Pour cela nous nous restreindrons au cas o\`u 
$k$ est un corps. 

Commen\c{c}ons par rappeler la d\'efinition d'une cat\'egorie
triangul\'ee $k$-lin\'eaire satur\'ee au sens de \cite[Def. 1.2]{bv}. 

\begin{df}\label{d11}
Une cat\'egorie triangul\'ee $k$-lin\'eaire $\mathcal{D}$ est \emph{satur\'ee}
si elle v\'erifie les deux conditions suivantes.
\begin{enumerate}
\item Pour toute paire d'objets $(x,y)$ dans $\mathcal{D}$, on a
$$Dim_{k}Ext^{*}(x,y)<\infty,$$
o\`u $Ext^{*}(x,y)=\oplus_{i\in \mathbb{Z}}[x,y[i]]$.
\item Tout foncteur cohomologique 
$$H : \mathcal{D}^{op} \longrightarrow Vect(k)$$
tel que pour tout $x\in \mathcal{D}$ on ait 
$$Dim_{k}\oplus_{i\in \mathbb{Z}}H(x[i]) < \infty,$$
est repr\'esentable. 
\end{enumerate}
\end{df}

Comme on peut s'y attendre la notion de saturation pour les dg-cat\'egories
est plus forte que celle pour les cat\'egories triangul\'ees. 

\begin{prop}\label{p13}
Soit $T$ une dg-cat\'egorie satur\'ee sur $k$. Alors la cat\'egorie
$[T]$, munie de sa structure triangul\'ee naturelle (voir \cite[\S 2.2]{tova}), est 
satur\'ee au sens de la d\'efinition \ref{d11}.
\end{prop}

\textit{Preuve:} Par d\'efinition on peut supposer que $T$ est de la forme
$L_{parf}(B)$, pour $B$ une dg-alg\`ebre propre et lisse sur $k$. Dans ce cas
$[T]\simeq D_{parf}(B)$ est la cat\'egorie triangul\'ee des $B$-dg-modules
parfaits (munie de sa structure triangul\'ee naturelle). On dispose d'un foncteur
bi-exact
$$D_{parf}(B\otimes B^{op}) \times D_{parf}(B) \longrightarrow D_{parf}(B),$$
qui envoie un $B\otimes B^{op}$-dg-module $P$ et un 
$B$-dg-module $M$ sur $P\otimes^{\mathbb{L}}_{B}M$. Ce foncteur
envoie $(B,M)$ sur $M$, et $(B\otimes B^{op},M)$ sur le 
$B$-dg-module libre $B\otimes M$. Comme $B$ appartient \`a l'enveloppe triangul\'ee
\'epaisse de $B\otimes B^{op}$, il existe un entier
$n$ tel que
$B\in <B\otimes B^{op}>_{n}$ (avec les notations de \cite{bv}). Cela implique facilement que
$M\in <B>_{n}$, et ce pour tout $M$. En d'autres termes, $B$ est un 
g\'en\'erateur fort (\emph{strong generator} dans la terminologie de \cite{bv}) 
de $D_{parf}(B)$. Comme $D_{parf}(B)$ satisfait la condition (1) de la
d\'efinition \ref{d11}, le th\'eor\`eme principal de \cite{bv} implique que 
$D_{parf}(B)$ est satur\'ee.
\hfill $\Box$ \\

Nous ne savons pas s'il est raisonnable d'attendre \`a ce qu'une r\'eciproque
de la proposition \ref{p13} soit vraie. On dispose cependant de la r\'eciproque partielle
suivante.

\begin{prop}\label{p14}
Soit $T$ une dg-cat\'egorie triangul\'ee (voir \cite[Def. 2.4]{tova}). On suppose que les
trois conditions suivantes sont satisfaites. 
\begin{enumerate}
\item La cat\'egorie triangul\'ee $[T]\simeq D_{parf}(T^{op})$ 
poss\`ede un g\'en\'erateur classique (au sens de \cite{bv}).
\item Pour toute paire d'objets $(x,y)$ dans $[T]$, on a
$$Dim_{k}Ext^{*}(x,y)<\infty,$$
o\`u $Ext^{*}(x,y)=\oplus_{i\in \mathbb{Z}}[x,y[i]]$.
\item La cat\'egorie triangul\'ee 
$[\widehat{T\otimes T^{op}}_{pe}]$ est satur\'ee. 
\end{enumerate}
Alors la dg-cat\'egorie $T$ est satur\'ee.
\end{prop}

\textit{Preuve:} La premi\`ere hypoth\`ese implique qu'il existe une
dg-alg\`ebre $B$ telle que $T$ soit quasi-\'equivalente 
\`a $L_{parf}(B)$. La seconde hypoth\`ese implique 
de plus que $B$ est une dg-alg\`ebre propre. 

La dg-cat\'egorie $\widehat{T\otimes T^{op}}_{pe}$ est 
quasi-\'equivalente \`a $L_{parf}(B\otimes B^{op})$. On consid\`ere
alors le foncteur cohomologique
$$H : D_{parf}(B\otimes B^{op}) \longrightarrow Vect(k),$$
qui \`a $M$ un $B\otimes B^{op}$dg-module parfait associe
$[M,B]$, l'ensemble des morphismes de $B$ vers $M$ pris
dans $D(B\otimes B^{op})$ la cat\'egorie d\'eriv\'ee 
de tous les $B\otimes B^{op}$-dg-modules. Ce foncteur 
v\'erifie bien la condition de finitude (2) de la d\'efinition \ref{d11} 
car $B\otimes B^{op}$ est un g\'en\'erateur classique
de $D_{parf}(B\otimes B^{op})$, et $H(B\otimes B^{op}[n])=H^{-n}(B)$
pour tout $n\in \mathbb{Z}$. Il existe donc un 
$B\otimes B^{op}$-dg-module $P$ parfait qui repr\'esente
le foncteur $H$. Choisissons un isomorphisme de foncteurs
$H\simeq [-,P]$. L'identit\'e de $P$ d\'efinit alors un morphisme 
dans $D(B\otimes B^{op})$
$$u : P \longrightarrow B.$$
Par construction, ce morphisme est tel que 
pour tout $M\in D_{parf}(B\otimes B^{op})$ le morphisme induit
$$u_{*} : [M,P] \longrightarrow [M,B]$$
soit bijectif. En appliquant cela aux $M=B\otimes B^{op}[n]$, on trouve
que pour tout $n\in \mathbb{Z}$ le morphisme induit
$$H^{n}(u) : H^{n}(P) \longrightarrow H^{n}(B)$$
est un isomorphisme. Ainsi, $u$ est un isomorphisme dans
$D(B\otimes B^{op})$ ce qui implique en particulier que $B$ 
est un objet de $D_{parf}(B\otimes B^{op})$. Par d\'efinition cela
signifie que $B$ est lisse. Ainsi, $T\simeq L_{parf}(B)$ avec
$B$ propre et lisse et donc $T$ est satur\'ee au sens de \cite{tova}.
\hfill $\Box$ \\

\section{Cat\'egories d\'eriv\'ees des espaces alg\'ebriques propres et lisses}

Le but de ce second appendice est de montrer la proposition suivante.

\begin{prop}\label{p15}
Soit $X$ un espace alg\'ebrique propre et lisse sur un corps $k$ de caract\'eristique nulle. Alors
la dg-cat\'egorie $L_{parf}(X)$ est  satur\'ee. 
\end{prop}

\textit{Preuve:} 
Lorsque $X$ est un sch\'ema cela est d\'emontr\'e dans
\cite[Lem. 3.27]{tova}. La m\^eme preuve marcherait aussi dans le cas g\'en\'eral si l'on 
savait que $D_{parf}(X)$ poss\`ede un g\'en\'erateur, ce que nous allons montrer.
Par le lemme de Chow et la r\'esolution des singularit\'es,  il existe un sch\'ema propre et lisse $X'$ et 
un morphisme birationnel surjectif $p : X' \longrightarrow X$. Ce morphisme
v\'erifie $\mathbb{R}p_{*}(\mathcal{O}_{X'})\simeq \mathcal{O}_{X}$ car 
$X$ est lisse et donc \`a singularit\'es rationnelles. Ceci implique par la formule de projection que
le foncteur
$$\mathbb{L}p^{*} : D_{parf}(X) \longrightarrow D_{parf}(X')$$
est pleinement fid\`ele. Mais cela implique en particulier que l'image d'un g\'en\'erateur
classique de $D_{parf}(X')$ par $\mathbb{R}p_{*}$ est un g\'en\'erateur classique de $D_{parf}(X)$. 
\hfill $\Box$ \\

\section{Cat\'egories d\'eriv\'ees non born\'ees}

Dans cette derni\`ere section nous d\'efinissons une cat\'egorie
de mod\`eles mono\"\i dale sym\'etrique dont la cat\'egorie
homotopique est la cat\'egorie d\'eriv\'ee non-born\'ee
des complexes de $\mathcal{O}$-modules sur un site
annel\'e $(C,\mathcal{O})$. L'existence de cette structure de mod\`eles
nous permet alors de consid\'erer les 4 op\'erations standards (images
directes, inverses, produits tensoriels et Hom internes) au niveau 
des cat\'egories d\'eriv\'ees non born\'ees.

Fixons $(C,\mathcal{O})$ un petit site annel\'e en anneaux commutatifs. 
Nous notons $C(\mathcal{O})$ la cat\'egorie des complexes
de pr\'efaisceaux de $\mathcal{O}$-modules sur $C$. Ses objets
sont les donn\'es de $\mathcal{O}(U)$-modules $M(U)$ pour tout objet
$U$ de $C$, fonctoriels en $U$ en un sens \'evident. Nous munissons
la cat\'egorie $C(\mathcal{O})$ d'une premi\`ere structure de mod\`eles, appel\'ee
la structure niveaux par niveaux pour la quelle on a:
\begin{enumerate}
\item Un morphisme $f : M \longrightarrow N$ est une
fibration si pour tout $U\in C$ le morphisme de complexes
de $\mathcal{O}(U)$-modules $M(U) \longrightarrow N(U)$ est un \'epimorphisme. 
\item Un morphisme $f : M \longrightarrow N$ est une
\'equivalence si pour tout $U\in C$ le morphisme
$M(U) \longrightarrow N(U)$ est un quasi-isomorphisme. 
\end{enumerate}

Il est facile de voir que ces notions d\'efinissent une structure de cat\'egorie de mod\`eles
engendr\'ee par cofibrations, propre et cellulaire sur $C(\mathcal{O})$. 

Nous introduisons maintenant la notion suivante d'\'equivalence locale. Pour cela rappelons
que pour $M\in C(\mathcal{O})$ et tout $n\in \mathbb{Z}$ on dispose de pr\'efaisceaux
$H^{n}_{pr}(M)$ sur $C$, qui \`a $U$ associe le $n$-\`eme groupe de
cohomologie du complexe $M(U)$. Le faisceau associ\'e \`a ce pr\'efaisceau sera 
not\'e $H^{n}(M)$. 

\begin{df}\label{d12}
Un morphisme $f : M \longrightarrow N$ dans $C(\mathcal{O})$ est une \emph{\'equivalence locale}
si pour tout $n\in \mathbb{Z}$ le morphisme induit
$$H^{n}(M) \longrightarrow H^{n}(N)$$
est un isomorphisme de faisceaux sur $C$.
\end{df}

La d\'efinition pr\'ec\'edente permet de d\'efinir des nouvelles notions
de cofibrations, fibrations et \'equivalences locales sur $C(\mathcal{O})$. 
Pour cela nous d\'efinissons les cofibrations locales comme \'etant
les cofibrations pour la structure niveaux par niveaux, et les fibrations locales
comme \'etant les morphismes poss\`edant la propri\'et\'e de rel\`evement
\`a droite des cofibrations qui sont aussi des \'equivalences locales. 

\begin{thm}\label{t3}
Les notions pr\'ec\'edente de fibrations, cofibrations et \'equivalences
locales d\'efinissent une structure de cat\'egorie de mod\`eles
sur $C(\mathcal{O})$, propre et engendr\'ee par cofibrations.
\end{thm}

Nous ne donnerons pas la preuve de ce th\'eor\`eme, elle est tout
\`a fait similaire \`a la preuve de l'existence de la structure locale
projective pour les pr\'efaisceaux simpliciaux donn\'ee par exemple
dans \cite[Thm. 3.4.1]{hagI}. \\

\begin{rmk}
\emph{La cat\'egorie homotopique $Ho(C(\mathcal{O}))$ est naturellement \'equivalente
\`a la cat\'egorie d\'eriv\'ee des faisceaux de $\mathcal{O}$-modules sur le site $C$. En effet, 
le foncteur d'inclusion des complexes de faisceaux de $\mathcal{O}$-modules dans
les complexes de pr\'efaisceaux de $\mathcal{O}$-modules induit clairement une \'equivalence
apr\`es avoir localis\'e le long des \'equivalences locales.}
\end{rmk}

La cat\'egorie $C(\mathcal{O})$ est naturellement munie d'une structure
mono\"\i dale sym\'etrique, not\'ee $\otimes$, et pour la quelle on a
$$(M\otimes N)(U)=M(U)\otimes N(U).$$

\begin{prop}\label{p16}
Supposons que $C$ poss\`ede des produits finis, alors
$C(\mathcal{O})$ munie de la structure mono\"\i dale 
$\otimes$ est une cat\'egorie de mod\`eles mono\"\i dale
au sens de \cite[\S 4]{ho}. Elle v\'erifie de plus l'axiome du mono\"\i de de \cite{ss}.
\end{prop}

\textit{Esquisse de preuve:} On commence par montrer que 
$C(\mathcal{O})$ est une cat\'egorie de mod\`eles mono\"\i dale
lorsqu'elle est munie de sa structure niveaux par niveaux. 
Pour cela, nous appliquerons l'\'enonc\'e g\'en\'eral de \cite[Cor. 4.2.5]{ho}, et 
nous aurons donc besoin d'expliciter les ensembles
g\'en\'erateurs $I$ et $J$ de cofibrations et cofibrations triviales.
Pour tout $U\in C$, le foncteur d'\'evaluation en $U$
$$j_{U}^{*} : C(\mathcal{O}) \longrightarrow C(\mathcal{O}(U)),$$
qui \`a $M$ associe le complexe de $\mathcal{O}(U)$-modules
$M(U)$, est de Quillen \`a droite (pour la structure
projective sur $C(\mathcal{O}(U))$ de \cite{ho}, pour la quelle les fibrations
sont les \'epimorphismes et les \'equivalences sont les quasi-isomorphismes). L'adjoint \`a gauche de ce foncteur
sera not\'e
$$(j_{U})_{!} : C(\mathcal{O}(U)) \longrightarrow C(\mathcal{O}).$$
Si $I(U)$ et $J(U)$ d\'esigne alors les ensembles g\'en\'erateurs de
cofibrations et cofibrations triviales dans $C(\mathcal{O}(U))$, 
il est facile de voir que 
$$I:=\{(j_{U})_{!}(u)\}_{U\in C,u\in I(U)} \qquad
J:=\{(j_{U})_{!}(u)\}_{U\in C,u\in J(U)}$$
sont des ensembles g\'en\'erateurs de cofibrations et cofibrations triviales 
pour la structure niveaux par niveaux de $C(\mathcal{O})$.

Soit maintenant 
$$(j_{U})_{!}(A) \longrightarrow (j_{U})_{!}(B) \qquad
(j_{V})_{!}(C) \longrightarrow (j_{V})_{!}(D)$$
deux \'el\'ements de $I$, pour deux objets $U,V\in C$. Alors, il est facile 
de voir qu'il existe un isomorphisme naturel
$$(j_{U})_{!}(X)\otimes (j_{V})_{!}(Y) \simeq (j_{U\times V})_{!}(X\otimes Y).$$
Ainsi, le morphisme induit
$$((j_{U})_{!}(A)\otimes (j_{V})_{!}(D)) \coprod_{
((j_{U})_{!}(A)\otimes (j_{V})_{!}(C))} ((j_{U})_{!}(B)\otimes (j_{V})_{!}(C)) \longrightarrow 
(j_{U})_{!}(B)\otimes (j_{V})_{!}(D)$$
est-il isomorphe au morphisme naturel
$$(j_{U\times V})_{!}((A\otimes D)\coprod_{A\otimes C}B\otimes C) \longrightarrow
(j_{U\times V})_{!}((B\otimes D)).$$
Ce dernier morphisme est bien une cofibration car $(j_{U\times V})_{!}$ est de Quillen 
\`a gauche. De plus cette cofibration est aussi une \'equivalence si l'un des deux morphismes
$A\longrightarrow B$ ou $C\longrightarrow D$ est une \'equivalence. Ceci fini de montrer 
que $C(\mathcal{O})$ munie de sa structure niveaux par niveaux est 
une cat\'egorie de mod\`eles mono\"\i dale. On en d\'eduit ais\'ement que
$C(\mathcal{O})$ munie de sa structure de mod\`eles locale est encore
une cat\'egorie de mod\`eles mono\"\i dale. 

Il reste \`a voir que $C(\mathcal{O})$ v\'erifie \`a l'axiome du mono\"\i de \cite{ss}. 
Pour cela il suffit de montrer que si $X \longrightarrow Y$ est une cofibration 
triviale locale dans $C(\mathcal{O})$, et $M$ est un objet de
$C(\mathcal{O})$ alors le morphisme 
$X\otimes M \longrightarrow Y\otimes M$ est une \'equivalence locale. On consid\`ere
la suite de morphismes dans $C(\mathcal{O})$
$$X\otimes M \longrightarrow Y\otimes M \longrightarrow Y/X\otimes M.$$
Comme $Y/X$ est cofibrant dans $C(\mathcal{O})$, il est plat niveaux par niveaux. Ainsi, 
la suite pr\'ec\'edente est une suite exacte, et il nous
suffit donc de montrer que $Y/X\otimes M$ est \'equivalent \`a $0$ (pour la structure
locale). Pour cela, soit $p : QM \longrightarrow M$ un remplacement cofibrant 
de $M$. Le morphisme $p$ est une fibration triviale et donc une \'equivalence
niveaux par niveaux. Comme $Y/X$ est plat niveaux par niveaux on voit que
le morphisme induit
$$Y/X\otimes QM \longrightarrow Y/X\otimes M$$
est une \'equivalence niveaux par niveaux. Ainsi, on a un isomorphisme dans $Ho(C(\mathcal{O}))$
$$Y/X \otimes^{\mathbb{L}}M \simeq Y/X \otimes M,$$
ce qui montre bien que $Y/X\otimes M$ est \'equivalent \`a $0$. 
\hfill $\Box$ \\

En corollaire de la proposition, on peut d\'efinir une cat\'egorie mono\"\i dale
sym\'etrique $D(C,\mathcal{O}):=Ho(C(\mathcal{O}))$, pour
la structure mono\"\i dale d\'eriv\'ee $\otimes^{\mathbb{L}}$. De plus, cette
structure mono\"\i dale est ferm\'ee. 
Le fait que l'axiome du mono\"\i de
soit v\'erifi\'e entraine l'existence de structures de mod\`eles
sur les cat\'egories de mono\"\i des et de modules dans $C(\mathcal{O})$ comme
d\'emontr\'e dans \cite{ss}.

Pour terminer, supposons que $f : C \longrightarrow D$ soit un 
foncteur exact \`a gauche et continu entre sites poss\'edant 
des limites finies. Supposons que $\mathcal{O}_{C}$ et 
$\mathcal{O}_{D}$ soient deux faisceaux d'anneaux commutatifs sur $C$ et $D$ munis d'un morphisme
$$\mathcal{O}_{C} \longrightarrow f^{-1}(\mathcal{O}_{D}).$$
On peut alors d\'efinir une adjonction
$$f_{!} : C(\mathcal{O}_{C}) \rightleftarrows C(\mathcal{O}_{D}) : f^{*},$$
dont il est facile de montrer qu'il s'agit d'une adjonction de Quillen en remarquant
que $f_{!}$ commute \`a la formation des pr\'efaisceaux de cohomologie 
$H_{pr}^{*}$.

\end{appendix}


\begin{thebibliography}{99}

\bibitem[An-To]{anto} M. Anel, B. To\"en, \textit{D\'enombrabilit\'e des
classes d'\'equivalences d\'eriv\'ees de vari\'et\'es alg\'ebriques}, 
Preprint math.AG/0611545.

\bibitem[Ar]{ar} M. Artin, \textit{
Algebraization of formal moduli II. Existence of modifications},
Ann. of Math. (2) $\mathbf{91}$ 1970 88--135. 

\bibitem[Ba-St]{bs} C. Banica, O. Stanasila, \textit{Algebraic methods in the 
global theory of complex spaces}. Translated from the Romanian. 
Editura Academiei, Bucharest; John Wiley  \& Sons, London-New York-Sydney, 1976. 296 pp.

\bibitem[B-V]{bv} A. Bondal, M. Van Den Bergh, \textit{Generators and representability of
functors in commutative and non-commutative geometry}, Mosc. Math. J. \textbf{3}
(2003), no. 1, 1--36.

\bibitem[Du]{du} D. Dugger, \textit{Universal Homotopy Theories}, Adv. in Math. \textbf{164} (2001), 144-176. 

\bibitem[D-H-I]{dhi}  D. Dugger, S. Hollander, D. Isaksen, \textit{Hypercovers and simplicial presheaves},
Math. Proc. Cambridge Philos. Soc.
\textbf{136} (2004), 9-51.

\bibitem[Gr-Re]{gr} H. Grauert, R.Remmert, 
\textit{Coherent analytic sheaves}, Grundlehren der Mathematischen Wissenschaften 
$\mathbf{265}$, Springer-Verlag, Berlin, 1984. xviii+249 pp.

\bibitem[Hi]{hi} P. Hirschhorn, \textit{Model categories and their localizations}, Math. Surveys and
Monographs Series \textbf{99}, AMS, Providence, 2003.

\bibitem[H-Si]{hs} A. Hirschowitz, C. Simpson, \textit{Descente pour les $n$-champs},
Preprint math.AG/9807049.

\bibitem[Ho]{ho} M. Hovey, \textit{Model categories}, Mathematical surveys and monographs, Vol. $\mathbf{63}$,
Amer. Math. Soc., Providence 1998.

\bibitem[Kn]{Kn} D. Knutson, \textit{Algebraic Spaces}, Lecture Notes in Mathematics $\mathbf{203}$, Springer-Verlag 1971.  

\bibitem[Ko-So]{koso} M. Kontsevich, Y. Soibelmann, \textit{Notes on A-infinity algebras, A-infinity 
categories and non-commutative geometry. I},
arXiv Preprint math.RA/0606241.

\bibitem[L]{l} M. Lieblich, \textit{Moduli of complexes on a proper morphism},
J. Algebraic Geom.  $\mathbf{15}$  (2006),  no. 1, 175--206.

\bibitem[SGA1]{sga1} \textit{Rev\^etements \'etales et groupe fondamental},  
[S\'eminaire de g\'eom\'etrie alg\'ebrique du Bois 
Marie 1960--61, dirig\'e par  
A. Grothendieck, Documents Math\'ematiques (Paris) 3, 
Soci\'et\'e Math\'ematique de France, Paris, 2003. xviii+327 pp.

\bibitem[S-S]{ss} S. Schwede, B. Shipley, \textit{Algebras and modules
in monoidal model categories}, Proc. London Math. Soc. (3)
$\mathbf{80}$ (2000), 491-511.

\bibitem[Si]{si} C. Simpson, \textit{Algebraic (geometric) $n$-stacks}, Preprint  math.AG/9609014.

\bibitem[To1]{to} B. To\"en, \textit{The homotopy theory of dg-caregories
and derived Morita theory}, \`a parraitre dans Invent. Math., 
pr\'e-publication math.AG/0408337.

\bibitem[To2]{to-dea} B. To\"en, \textit{Champs alg\'ebriques: cours 2}, cours de Master 2 accessible
\`a http://www.picard.ups-tlse.fr/~toen/m2.html.

\bibitem[To-Va]{tova} B. To\"en, M. Vaqui\'e, \textit{Moduli of objects in dg-categories}, 
pr\'e-publication arXiv math.AG/0503269.

\bibitem[HAGI]{hagI} B. To\"en, G. Vezzosi, \textit{Homotopical algebraic geometry I: Topos theory}, 
Adv. in Math. \textbf{193}, (2005), 257-372.  

\bibitem[HAGII]{hagII} B. To\"en, G. Vezzosi,
\textit{Homotopical algebraic geometry II: Geometric stacks and applications},
\`a parraitre dans M\'emoires of the AMS,
pr\'e-publication math.AG/0404373.



\end{thebibliography}
\end{document}